# ON THE EXISTENCE OF WEAK SOLUTIONS OF SEMILINEAR ELLIPTIC EQUATIONS AND SYSTEMS WITH HARDY POTENTIALS

KONSTANTINOS T. GKIKAS AND PHUOC-TAI NGUYEN

Abstract. Let $\Omega \subset \mathbb{R}^N$ $(N \geq 3)$ be a bounded $C^2$ domain and $\delta(x) = \mathrm{dist}\,(x, \partial\Omega)$. Put $L_\mu = \Delta + \frac{\mu}{\delta^2}$ with $\mu > 0$. In this paper, we provide various necessary and sufficient conditions for the existence of weak solutions to

$$-L_\mu u = u^p + \tau \quad \text{in } \Omega, \qquad u = \nu \quad \text{on } \partial\Omega,$$

where $\mu > 0$, $p > 0$, $\tau$ and $\nu$ are measures on $\Omega$ and $\partial\Omega$ respectively. We then establish existence results for the system

$$\begin{cases} -L_\mu u = \epsilon\, v^p + \tau & \text{in } \Omega, \\ -L_\mu v = \epsilon\, u^{\tilde{p}} + \tilde{\tau} & \text{in } \Omega, \\ u = \nu, \quad v = \tilde{\nu} & \text{on } \partial\Omega, \end{cases}$$

where $\epsilon = \pm 1$, $p > 0$, $\tilde{p} > 0$, $\tau$ and $\tilde{\tau}$ are measures on $\Omega$, $\nu$ and $\tilde{\nu}$ are measures on $\partial\Omega$. We also deal with elliptic systems where the nonlinearities are more general.

## Contents

Email address: kgkikas@dim.uchile.cl.

Email address: ptnguyen@math.muni.cz; nguyenphuoctai.hcmup@gmail.com.

## 1. Introduction

Let $\Omega \subset \mathbb{R}^N$ $(N \geq 3)$ be a bounded $C^2$ domain, $\delta(x) = \mathrm{dist}\,(x, \partial\Omega)$ and $g \in C(\mathbb{R})$. Put $L_\mu := \Delta + \frac{\mu}{\delta^2}$. In the present paper we study semilinear problems with *Hardy potential* of the form

$$-L_\mu u = g(u) + \tau \quad \text{in } \Omega, \tag{1.1}$$

where $\mu > 0$, $\tau$ is a Radon measure on $\Omega$.

The boundary value problem with measures for (1.1) without Hardy potential and with power absorption nonlinearity, i.e. $\mu = 0$, $\tau = 0$, $g(u) = -|u|^{p-1}u$, $p > 1$, was well understood in the literature, starting with a work by Gmira and Véron [10]. It was proved that there is the critical exponent $p^* := \frac{N+1}{N-1}$ in the sense that if $p \in (1, p^*)$ then there is a unique weak solution for every finite measure $\nu$ on $\Omega$, while if $p \in [p^*, \infty)$ there exists no solution with a boundary isolated singularity. Marcus and Véron [15] studied this problem by introducing a notion of *boundary trace*, providing a complete description of isolated singularities in the *subcritical case*, i.e. $1 < p < p^*$, and giving a removability result in the *supercritical case*, i.e. $p \geq p^*$.

The solvability for boundary value problem for (1.1) without Hardy potential and with *power source term*, namely $\mu = 0$, $\tau = 0$, $g(u) = u^p$, $p > 1$, was studied by Bidaut-Véron and Vivier [4] in connection with sharp estimates of the Green operator and the Poisson operator associated to $(-\Delta)$ in $\Omega$. They proved that, in the subcritical case $1 < p < p^*$, the problem admits a solution if and only if the total mass of the boundary datum $\nu$ is sufficiently small. Afterwards, Bidaut-Véron and Yarur [6] reconsidered this type of problem in a more general setting and provided a necessary and sufficient condition for the existence of solutions. Recently, Bidaut-Véron et al. [5] provided new criteria for the existence of solutions with $p > 1$ in terms of the capacity associated to the Besov spaces.

Let $\phi \geq 0$ in $\Omega$ and $p \geq 1$, we denote by $L^p(\Omega; \phi)$ the space of all function $v$ on $\Omega$ satisfying $\int_\Omega |v|^p \, \phi dx < \infty$. We denote by $\mathfrak{M}(\Omega; \phi)$ the space of Radon measures $\tau$ on $\Omega$ satisfying $\int_\Omega \phi \, d|\tau| < \infty$ and by $\mathfrak{M}^+(\Omega; \phi)$ the nonnegative cone of $\mathfrak{M}(\Omega; \phi)$. When $\phi \equiv 1$, we use the notations $\mathfrak{M}(\Omega)$ and $\mathfrak{M}^+(\Omega)$. We also denote by $\mathfrak{M}(\partial\Omega)$ the space of finite measures on $\partial\Omega$ and by $\mathfrak{M}^+(\partial\Omega)$ the nonnegative cone of $\mathfrak{M}(\partial\Omega)$.

Let $G_\mu$ and $K_\mu$ be the Green kernel and Martin kernel of $-L_\mu$ in $\Omega$, $\mathbb{G}_\mu$ and $\mathbb{K}_\mu$ be the corresponding Green operator and Martin operator (see [14, 9] for more details). Let $C_H$ be Hardy constant, namely

$$C_H := \inf_{v \in H_0^1(\Omega)\setminus\{0\}} \frac{\int_\Omega |\nabla v|^2 dx}{\int_\Omega (v/\delta)^2 dx} \tag{1.2}$$

then it is well known that $0 < C_H \leq \frac{1}{4}$ and if $\Omega$ is convex then $C_H = \frac{1}{4}$ (see for example [12]). Moreover the infimum is achieved if and only if $C_H < \frac{1}{4}$. When $-\Delta\delta \geq 0$ in $\Omega$ in the sense of distributions, the first eigenvalue $\lambda_\mu$ of $L_\mu$ in $\Omega$ is positive, i.e.

$$\lambda_\mu := \inf_{\varphi \in H_0^1(\Omega)\setminus\{0\}} \frac{\int_\Omega (|\nabla \varphi|^2 - \frac{\mu}{\delta^2}\varphi^2) dx}{\int_\Omega \varphi^2 dx} > 0. \tag{1.3}$$

For $\mu \in (0, \frac{1}{4}]$, denote by $\alpha$ the following *fundamental exponent*

$$\alpha := \frac{1}{2}(1 + \sqrt{1 - 4\mu}). \tag{1.4}$$

Notice that $\frac{1}{2} < \alpha < 1$. The eigenfunction $\varphi_\mu$ associated to $\lambda_\mu$ with the normalization $\int_\Omega (\varphi_\mu/\delta)^2 dx = 1$ satisfies $c^{-1}\delta^\alpha \le \varphi_\mu \le c\delta^\alpha$ for some constant $c > 0$ (see [7]).

In relation to Hardy constant, Bandle et al. [3] classified large solutions of the linear equation

$$-L_\mu u = 0 \quad \text{in } \Omega, \tag{1.5}$$

and of the associated nonlinear equation with power absorption

$$-L_\mu u + u^p = 0 \quad \text{in } \Omega. \tag{1.6}$$

In [14], Marcus and P.-T. Nguyen studied boundary value problem for (1.5) and (1.6) with $\mu \in (0, C_H)$ in measure framework by introducing a notion of *normalized boundary trace* which is defined as follows:

**Definition 1.1.** *A function* $u \in L^1_{loc}(\Omega)$ *possesses a* normalized boundary trace *if there exists a measure* $\nu \in \mathfrak{M}(\partial\Omega)$ *such that*

$$\lim_{\beta \to 0} \beta^{\alpha-1} \int_{\{x \in \Omega : \delta(x) = \beta\}} |u - \mathbb{K}_\mu[\nu]| dS = 0. \tag{1.7}$$

*The normalized boundary trace is denoted by* $\operatorname{tr}^*(u)$*.*

The restriction $\mu \in (0, C_H)$ in [14] is due to the fact that in this case $L_\mu$ is *weakly coercive* in $H^1_0(\Omega)$ and consequently by a result of Ancona [2, Remark p. 523] there is a $(1-1)$ correspondence between $\mathfrak{M}^+(\partial\Omega)$ and the class of positive $L_\mu$ harmonic functions, namely any positive $L_\mu$ harmonic function $u$ can be written in a unique way under the form $u = \mathbb{K}_\mu[\nu]$ for some $\nu \in \mathfrak{M}^+(\partial\Omega)$.

The notion of normalized boundary trace was proved [14] to be an appropriate generalization of the classical boundary trace to the setting of Hardy potentials, giving a characterization of *moderate solutions* of (1.6). In addition, it was showed in [14] that there exists the critical exponent

$$p_\mu := \frac{N + \alpha}{N + \alpha - 2} \tag{1.8}$$

such that if $p \in (1, p_\mu)$ then there exists a unique solution of (1.6) with $\operatorname{tr}^*(u) = \nu$ for every finite measure $\nu$ on $\partial\Omega$, while if $p \ge p_\mu$ there is no solution of (1.6) with an isolated boundary singularity. Marcus and Moroz [13] then *extended the notion of normalized boundary trace to the case* $\mu < \frac{1}{4}$ and employed it to investigate (1.6). When $\mu = \frac{1}{4}$, $L_\mu$ is no longer weakly coercive and hence Ancona's result cannot be applied. However, Gkikas and Véron [9] observed that *if the first eigenvalue of* $-L_{\frac{1}{4}}$ *is positive* then the kernel $K_{\frac{1}{4}}(\cdot, y)$ with pole at $y \in \partial\Omega$ is unique up to a multiplication and any positive $L_{\frac{1}{4}}$ harmonic function $u$ admits such a representation. Based on that observation, they considered the boundary value problem with measures for (1.6), fully classifying isolated singularities in the subcritical case $p \in (1, p_\mu)$ and providing removability result in the supercritical case $p \ge p_\mu$. A main ingredient in [9] is the notion of boundary trace which is defined in a *dynamic way* and is recalled below.

Let $D \Subset \Omega$ and $x_0 \in D$. If $h \in C(\partial D)$ then the following problem

$$\begin{cases} -L_\mu u = 0 & \text{in } D, \\ u = h & \text{on } \partial D, \end{cases} \tag{1.9}$$

admits a unique solution which allows to define the $L_\mu$-harmonic measure $\omega_D^{x_0}$ on $\partial D$ by

$$u(x_0) = \int_{\partial D} h(y) d\omega_D^{x_0}(y). \tag{1.10}$$

A sequence of domains $\{\Omega_n\}$ is called a *smooth exhaustion* of $\Omega$ if $\partial\Omega_n \in C^2$, $\overline{\Omega_n} \subset \Omega_{n+1}$, $\cup_n \Omega_n = \Omega$ and $\mathcal{H}^{N-1}(\partial\Omega_n) \to \mathcal{H}^{N-1}(\partial\Omega)$. For each $n$, let $\omega_{\Omega_n}^{x_0}$ be the $L_\mu^{\Omega_n}$-harmonic measure on $\partial\Omega_n$.

**Definition 1.2.** *A function $u$ possesses a* boundary trace *if there exists a measure $\nu \in \mathfrak{M}(\partial\Omega)$ such that for any smooth exhaustion $\{\Omega_n\}$ of $\Omega$,*

$$\lim_{n\to\infty} \int_{\partial\Omega_n} \zeta u \, d\omega_{\Omega_n}^{x_0} = \int_{\partial\Omega} \zeta \, d\nu \quad \forall \zeta \in C(\overline{\Omega}). \tag{1.11}$$

*The boundary trace of $u$ is denoted by* $\mathrm{tr}\,(u)$.

It is worthy mentioning that in Definition 1.2, $\mu$ is allowed to belong to the range $(0, \frac{1}{4}]$ provided $\lambda_\mu > 0$.

In parallel, semilinear equations with Hardy potential and source term

$$-L_\mu u = u^p \quad \text{in } \Omega \tag{1.12}$$

were treated by Bidaut-Véron et al. [5] and by P.-T. Nguyen [18] and a fairly complete description of the profile of solutions to (1.12) was obtained in subcritical case $p < p_\mu$ (see [18]) and in supercritical case $p \geq p_\mu$ (see [5]).

Our first contribution is to show that *the notion of normalized boundary trace given in Definition 1.1 is equivalent to that in Definition 1.2 by examining* $\mathrm{tr}\,(\mathbb{G}_\mu^\Omega[\tau]) = \mathrm{tr}^*(\mathbb{G}_\mu^\Omega[\tau])$ *and* $\mathrm{tr}\,(\mathbb{K}_\mu^\Omega[\nu]) = \mathrm{tr}^*(\mathbb{K}_\mu^\Omega[\nu])$. This enables to establish important results for the boundary value problem for linear equations (see Proposition 2.13) which in turn forms a basic to study the boundary value problem for

$$\begin{cases} -L_\mu u = g(u) + \tau & \text{in } \Omega, \\ \mathrm{tr}\,(u) = \nu. \end{cases} \tag{1.13}$$

When dealing with (1.13), one encounters the following difficulties. The first one is due to the presence of the Hardy potential in the linear part of the equations. More precisely, since the singularity of the potential at the boundary is too strong, some important tools such as Hopf's lemma, the classical notion of boundary trace, are invalid, and therefore the system cannot be handled via classical elliptic PDEs methods. The second one comes from the interplay between the nonlinearity, the Hardy potential and measure data. The interaction between the difficulties generates an intricate dynamics both in $\Omega$ and near $\partial\Omega$ and leads to disclose new type of results.

**Convention.** Throughout the paper, unless otherwise stated, we assume that $\mu \in (0, \frac{1}{4}]$ *and the first eigenvalue $\lambda_\mu$ of $-L_\mu$ in $\Omega$ is positive.* We emphasize that if $\mu \in (0, C_H)$ then $\lambda_\mu > 0$.

**Definition 1.3.** *(i) The space of test functions is defined as*

$$\mathbf{X}_\mu(\Omega) := \{\zeta \in H^1_{loc}(\Omega) : \delta^{-\alpha}\zeta \in H^1(\Omega, \delta^{2\alpha}),\, \delta^{-\alpha} L_\mu \zeta \in L^\infty(\Omega)\}. \tag{1.14}$$

*(ii) Let $(\tau,\nu) \in \mathfrak{M}(\Omega,\delta^\alpha) \times \mathfrak{M}(\partial\Omega)$. We say that $u$ is a* weak solution *of* (1.13) *if $u \in L^1(\Omega;\delta^\alpha)$, $g(u) \in L^1(\Omega;\delta^\alpha)$ and*

$$-\int_\Omega u L_\mu \zeta \, dx = \int_\Omega g(u)\zeta \, dx + \int_\Omega \zeta d\tau - \int_\Omega \mathbb{K}_\mu[\nu] L_\mu \zeta \, dx \quad \forall \zeta \in \mathbf{X}_\mu(\Omega). \tag{1.15}$$

Main properties of solutions of (1.13) are established in the following proposition.

**Proposition A.** *Let $\tau \in \mathfrak{M}(\Omega;\delta^\alpha)$ and $\nu \in \mathfrak{M}(\partial\Omega)$. The following statements are equivalent.*
(i) *$u$ is a weak solution of* (1.13).
(ii) *$g(u) \in L^1(\Omega;\delta^\alpha)$ and*

$$u = \mathbb{G}_\mu[g(u)] + \mathbb{G}_\mu[\tau] + \mathbb{K}_\mu[\nu]. \tag{1.16}$$

(iii) *$u \in L^1_{loc}(\Omega)$, $g(u) \in L^1_{loc}(\Omega)$, $u$ is a distributional solution of* (1.1) *and* $\mathrm{tr}\,(u) = \nu$.

This allows to establish necessary and sufficient conditions for the existence of a weak solution of

$$\begin{cases} -L_\mu u = u^p + \sigma\tau & \text{in } \Omega, \\ \mathrm{tr}\,(u) = \varrho\nu. \end{cases} \tag{1.17}$$

**Theorem B.** *Let $\tau \in \mathfrak{M}^+(\Omega;\delta^\alpha)$, $\nu \in \mathfrak{M}^+(\partial\Omega)$ and $p > 0$.*
(i) *Assume $0 < p < p_\mu$. Then there exists a constant $C > 0$ such that*

$$\mathbb{G}_\mu[\mathbb{K}_\mu[\nu]^p] \le C\mathbb{K}_\mu[\nu] \quad a.e. \text{ in } \Omega. \tag{1.18}$$

(ii) *Assume $0 < p < p_\mu$. Then there exists a constant $C > 0$ such that*

$$\mathbb{G}_\mu[\mathbb{G}_\mu[\tau]^p] \le C\mathbb{G}_\mu[\tau] \quad a.e. \text{ in } \Omega. \tag{1.19}$$

(iii) *If* (1.18) *and* (1.19) *hold then problem* (1.17) *admits a weak solution $u$ satisfying*

$$\mathbb{G}_\mu[\sigma\tau] + \mathbb{K}_\mu[\varrho\nu] \le u \le C(\mathbb{G}_\mu[\sigma\tau] + \mathbb{K}_\mu[\varrho\nu]) \ a.e. \text{ in } \Omega \tag{1.20}$$

*for $\sigma > 0$ and $\varrho > 0$ small enough if $p > 1$, for any $\sigma > 0$ and $\varrho > 0$ if $0 < p < 1$.*
(iv) *If $p > 1$ and* (1.17) *admits a weak solution then* (1.18) *and* (1.19) *hold with constant $C = \frac{1}{p-1}$.*
(v) *Assume $0 < p < p_\mu$. Then there exists a constant $C > 0$ such that for any weak solution $u$ of* (1.17) *there holds*

$$\mathbb{G}_\mu[\sigma\tau] + \mathbb{K}_\mu[\varrho\nu] \le u \le C(\mathbb{G}_\mu[\sigma\tau] + \mathbb{K}_\mu[\varrho\nu] + \delta^\alpha) \ a.e. \text{ in } \Omega. \tag{1.21}$$

In order to study (1.17) in the supercritical case, i.e. $p \ge p_\mu$, we make use of the capacities introduced in [5] which is recalled below. For $0 \le \theta \le \beta < N$, set

$$N_{\theta,\beta}(x,y) := \frac{1}{|x-y|^{N-\beta}\max\{|x-y|,\delta(x),\delta(y)\}^\theta} \quad \forall (x,y) \in \overline{\Omega} \times \overline{\Omega},\ x \ne y, \tag{1.22}$$

$$\mathbb{N}_{\theta,\beta}[\tau](x) := \int_{\overline{\Omega}} N_{\theta,\beta}(x,y) d\tau \quad \forall \tau \in \mathfrak{M}^+(\overline{\Omega}). \tag{1.23}$$

For $a > -1$, $0 \le \theta \le \beta < N$ and $s > 1$, define $\mathrm{Cap}^a_{\mathbb{N}_{\theta,\beta},s}$ by

$$\mathrm{Cap}^a_{\mathbb{N}_{\theta,\beta},s}(E) := \inf\left\{ \int_{\overline{\Omega}} \delta^a \phi^s \, dx : \quad \phi \ge 0,\ \ \mathbb{N}_{\theta,\beta}[\delta^a\phi] \ge \chi_E \right\}, \tag{1.24}$$

for any Borel set $E \subset \overline{\Omega}$. For $\theta \in (0, N-1)$ and $s > 0$, let $\mathrm{Cap}^{\partial\Omega}_{\theta,s}$ be the capacity defined in [5, Definition 1.1]. Notice that if $\theta s > N-1$ then $\mathrm{Cap}^{\partial\Omega}_{\theta,s}(\{z\}) > 0$ for every $z \in \partial\Omega$.

**Theorem C.** *Let $\tau \in \mathfrak{M}^+(\Omega;\delta^\alpha)$ and $\nu \in \mathfrak{M}^+(\partial\Omega)$. Assume $p > 1$. Then the following statements are equivalent.*

(i) *There exists $C > 0$ such that the following inequalities hold*

$$\int_E \delta^\alpha d\tau \leq C\mathrm{Cap}^{(p+1)\alpha}_{\mathbb{N}_{2\alpha,2,p'}}(E) \quad \forall \text{ Borel } E \subset \overline{\Omega}, \tag{1.25}$$

$$\nu(F) \leq C\mathrm{Cap}^{\partial\Omega}_{1-\alpha+\frac{\alpha+1}{p},p'}(F) \quad \forall \text{ Borel } F \subset \partial\Omega. \tag{1.26}$$

(ii) *There exists a positive constant $C$ such that* (1.18) *and* (1.19) *hold.*

(iii) *Problem* (1.17) *has a positive weak solution for $\sigma > 0$ and $\varrho > 0$ small enough.*

**Remark.** When $\tau = 0$, Theorem C covers Theorem B (i), (iii) due to the fact that $\mathrm{Cap}^{\partial\Omega}_{1-\alpha+\frac{\alpha+1}{p},p'}(\{z\}) > c > 0$ for every $z \in \partial\Omega$ if $1 < p < p_\mu$. Also if $1 < p < p_\mu$ then (see Lemma 3.10)

$$\inf_{\xi\in\Omega} \mathrm{Cap}^{(p+1)\alpha}_{\mathbb{N}_{2\alpha,2,p'}}(\{\xi\}) > 0,$$

which implies the statements (ii) and (iii) in Theorem B.

The next goal of the present paper is the study of weak solutions of semilinear elliptic system involving Hardy potential

$$\begin{cases} -L_\mu u = g(v) + \tau & \text{in } \Omega, \\ -L_\mu v = \tilde{g}(u) + \tilde{\tau} & \text{in } \Omega, \\ \mathrm{tr}\,(u) = \nu, \quad \mathrm{tr}\,(v) = \tilde{\nu} \end{cases} \tag{1.27}$$

where $\tau,\ \tilde{\tau} \in \mathfrak{M}(\Omega;\delta^\alpha)$, $\nu, \tilde{\nu} \in \mathfrak{M}(\partial\Omega)$, $g, \tilde{g} \in C(\mathbb{R})$.

**Definition 1.4.** *A pair $(u, v)$ is called a weak solution of* (1.27) *if $u \in L^1(\Omega;\delta^\alpha)$, $v \in L^1(\Omega;\delta^\alpha)$, $\tilde{g}(u) \in L^1(\Omega;\delta^\alpha)$, $g(v) \in L^1(\Omega;\delta^\alpha)$ and*

$$\begin{aligned} -\int_\Omega uL_\mu\zeta\, dx &= \int_\Omega g(v)\zeta\, dx + \int_\Omega \zeta\, d\tau - \int_\Omega \mathbb{K}_\mu[\nu]L_\mu\zeta\, dx, \\ -\int_\Omega vL_\mu\zeta\, dx &= \int_\Omega \tilde{g}(u)\zeta\, dx + \int_\Omega \zeta\, d\tilde{\tau} - \int_\Omega \mathbb{K}_\mu[\tilde{\nu}]L_\mu\zeta\, dx \quad \forall \zeta \in \mathbf{X}_\mu(\Omega). \end{aligned} \tag{1.28}$$

A counterpart of Proposition A in the case of systems is the following:

**Proposition D.** *Let $\tau, \tilde{\tau} \in \mathfrak{M}(\Omega;\delta^\alpha)$ and $\nu, \tilde{\nu} \in \mathfrak{M}(\partial\Omega)$. Then the following statements are equivalent.*

(i) *$(u, v)$ is a weak solution of* (1.27).

(ii) *$\tilde{g}(u) \in L^1(\Omega;\delta^\alpha)$, $g(v) \in L^1(\Omega;\delta^\alpha)$ and*

$$u = \mathbb{G}_\mu[g(v)] + \mathbb{G}_\mu[\tau] + \mathbb{K}_\mu[\nu], \quad v = \mathbb{G}_\mu[\tilde{g}(u)] + \mathbb{G}_\mu[\tilde{\tau}] + \mathbb{K}_\mu[\tilde{\nu}]. \tag{1.29}$$

(iii) *$(u, v) \in (L^1_{loc}(\Omega))^2$, $(g(v), \tilde{g}(u)) \in (L^1_{loc}(\Omega))^2$, $(u, v)$ is a solution of*

$$\begin{cases} -L_\mu u = g(v) + \tau & \text{in } \Omega, \\ -L_\mu v = \tilde{g}(u) + \tilde{\tau} & \text{in } \Omega, \end{cases} \tag{1.30}$$

*in the sense of distributions and* $\mathrm{tr}\,(u) = \nu$ *and* $\mathrm{tr}\,(v) = \tilde{\nu}$.

Elliptic systems arise in biological applications (e.g. population dynamics) or physical applications (e.g. models of nuclear reactor) and have been drawn a lot of attention (see [8, 19] and references therein). A typical case is Lane-Emden system, i.e. system (1.27) with $\mu = 0$, $g(v) = v^p$, $\tilde{g}(u) = u^{\tilde{p}}$. Bidaut-Véron and Yarur [6] proved various existence results for Lane-Emden system under conditions involving the following exponents

$$q := \tilde{p}\frac{p+1}{\tilde{p}+1}, \quad \tilde{q} := p\frac{\tilde{p}+1}{p+1}. \tag{1.31}$$

We first treat the system

$$\begin{cases} -L_\mu u = v^p + \sigma\tau & \text{in } \Omega, \\ -L_\mu v = \tilde{u}^{\tilde{p}} + \tilde{\sigma}\tilde{\tau} & \text{in } \Omega, \\ \mathrm{tr}\,(u) = \varrho\nu, \quad \mathrm{tr}\,(v) = \tilde{\varrho}\tilde{\nu}, \end{cases} \tag{1.32}$$

where $p > 0, \tilde{p} > 0$, $\tau, \tilde{\tau} \in \mathfrak{M}(\Omega; \delta^\alpha)$ and $\nu, \tilde{\nu} \in \mathfrak{M}(\partial\Omega)$.

The next theorem provides a sufficient condition for the existence of solutions of (1.32).

**Theorem E.** *Let $p > 0$, $\tilde{p} > 0$, $\tau, \tilde{\tau} \in \mathfrak{M}^+(\Omega; \delta^\alpha)$ and $\nu, \tilde{\nu} \in \mathfrak{M}^+(\partial\Omega)$. Assume $p\tilde{p} \neq 1$, $q < p_\mu$, $\mathbb{G}_\mu[\tau] + \mathbb{K}_\mu[\nu + \tilde{\nu}] \in L^{\tilde{p}}(\Omega, \delta^\alpha)$. Then system (1.32) admits a weak solution $(u, v)$ for $\sigma > 0$ and $\tilde{\sigma} > 0$ small if $p\tilde{p} > 1$, for any $\sigma > 0$ and $\tilde{\sigma} > 0$ if $p\tilde{p} < 1$. Moreover*

$$v \approx \mathbb{G}_\mu[\omega] + \mathbb{K}_\mu[\tilde{\nu}], \tag{1.33}$$

$$u \approx \mathbb{G}_\mu[(\mathbb{G}_\mu[\omega] + \mathbb{K}_\mu[\tilde{\nu}])^p] + \mathbb{G}_\mu[\tau] + \mathbb{K}_\mu[\nu] \tag{1.34}$$

*where the similarity constants depend on $N, p, \tilde{p}, \mu, \Omega, \sigma, \tilde{\sigma}, \tau, \tilde{\tau}$ and*

$$\omega := \mathbb{G}_\mu[\tau + \mathbb{K}_\mu[\tilde{\nu}]^p]^{\tilde{p}} + \mathbb{K}_\mu[\nu]^{\tilde{p}} + \tilde{\tau}.$$

A new criterion for the existence of (1.32), expressed in terms of the capacities $\mathrm{Cap}^a_{\mathbb{N}_{\theta,\beta},s}$ and $\mathrm{Cap}^{\partial\Omega}_{\theta,s}$, is stated in the following result.

**Theorem F.** *Let $p > 1$, $\tilde{p} > 1$, $\tau, \tilde{\tau} \in \mathfrak{M}^+(\Omega; \delta^\alpha)$ and $\nu, \tilde{\nu} \in \mathfrak{M}^+(\partial\Omega)$. Assume there exists $C > 0$ such that*

$$\max\left\{\int_E \delta^\alpha d\tau, \int_E \delta^\alpha d\tilde{\tau}\right\} \leq C \min\{\mathrm{Cap}^{(p+1)\alpha}_{\mathbb{N}_{2\alpha,2},p'}(E), \mathrm{Cap}^{(\tilde{p}+1)\alpha}_{\mathbb{N}_{2\alpha,2},\tilde{p}'}(E)\}, \ \forall E \subset \overline{\Omega}, \tag{1.35}$$

$$\max\{\nu(F), \tilde{\nu}(F)\} \leq C \min\{\mathrm{Cap}^{\partial\Omega}_{1-\alpha+\frac{1+\alpha}{p},p'}(F), \mathrm{Cap}^{\partial\Omega}_{1-\alpha+\frac{1+\alpha}{\tilde{p}},\tilde{p}'}(F)\}, \ \forall F \subset \partial\Omega. \tag{1.36}$$

*Then (1.32) admits a weak solution $(u, v)$ for $\sigma > 0$, $\tilde{\sigma} > 0$, $\varrho > 0$, $\tilde{\varrho} > 0$ small enough. There exists $C > 0$ such that*

$$\begin{aligned} \mathbb{G}_\mu[\sigma\tau] + \mathbb{K}_\mu[\varrho\nu] &\leq u \leq C(\mathbb{G}_\mu[\sigma\tau + \tilde{\sigma}\tilde{\tau}] + \mathbb{K}_\mu[\varrho\nu + \tilde{\varrho}\tilde{\nu}]), \\ \mathbb{G}_\mu[\tilde{\sigma}\tilde{\tau}] + \mathbb{K}_\mu[\tilde{\varrho}\tilde{\nu}] &\leq v \leq C(\mathbb{G}_\mu[\sigma\tau + \tilde{\sigma}\tilde{\tau}] + \mathbb{K}_\mu[\varrho\nu + \tilde{\varrho}\tilde{\nu}]). \end{aligned} \tag{1.37}$$

Finally, we deal with elliptic systems with more general nonlinearities

$$\begin{cases} -L_\mu u = \epsilon\, g(v) + \sigma\tau & \text{in } \Omega, \\ -L_\mu v = \epsilon\, \tilde{g}(u) + \tilde{\sigma}\tilde{\tau} & \text{in } \Omega, \\ \mathrm{tr}\,(u) = \varrho\nu, \quad \mathrm{tr}\,(v) = \tilde{\varrho}\tilde{\nu} & \text{on } \partial\Omega \end{cases} \tag{1.38}$$

where $g$ and $\tilde{g}$ are nondecreasing, continuous functions in $\mathbb{R}$, $\epsilon = \pm 1$, $\sigma > 0$, $\tilde{\sigma} > 0$, $\varrho > 0$, $\tilde{\varrho} > 0$.

We shall treat successively the cases $\epsilon = -1$ and $\epsilon = 1$. For any function $f$, define

$$\Lambda_f := \int_1^\infty s^{-1-p_\mu}|f(s) - f(-s)|ds \tag{1.39}$$

with $p_\mu$ defined in (1.8).

**Theorem G.** *Let $\epsilon = -1$ and $\sigma, \tilde\sigma, \varrho, \tilde\varrho$ be positive numbers, $\tau,\ \tilde\tau \in \mathfrak{M}(\Omega;\delta^\alpha)$ and $\nu, \tilde\nu \in \mathfrak{M}(\partial\Omega)$. Assume that $\Lambda_g + \Lambda_{\tilde g} < \infty$ and $g(s) = \tilde g(s) = 0$ for any $s \le 0$. Then system* (1.38) *admits a weak solution $(u, v)$.*

When $\epsilon = 1$, different phenomenon occurs, which is reflected in the following result.

**Theorem H.** *Let $\epsilon = 1$, $\tau,\ \tilde\tau \in \mathfrak{M}(\Omega;\delta^\alpha)$ and $\nu, \tilde\nu \in \mathfrak{M}(\partial\Omega)$.*
I. Subcriticality. *Assume that $\Lambda_g + \Lambda_{\tilde g} < \infty$. In addition, assume that there exist $q_1 > 1$, $a_1 > 0$, $b_1 > 0$ such that*

$$|g(s)| \le a_1|s|^{q_1} + b_1 \quad \forall s \in [-1, 1], \tag{1.40}$$

$$|\tilde g(s)| \le a_1|s|^{q_1} + b_1 \quad \forall s \in [-1, 1]. \tag{1.41}$$

*Then* (1.38) *admits a weak solution for $b_1, \sigma, \tilde\sigma, \varrho, \tilde\varrho$ small enough.*

II. Sublinearity. *Assume that there exist $q_1 > 1$, $q_2 \in (0, 1]$, $a_2 > 0$ and $b_2 > 0$ such that $\mathbb{K}_\mu[|\tilde\nu|] + \mathbb{G}_\mu[|\tilde\tau|] \in L^{q_1}(\Omega;\delta^{\alpha-1})$ and*

$$|g(s)| \le a_2|s|^{q_1} + b_2 \quad \forall s \in \mathbb{R}, \tag{1.42}$$

$$|\tilde g(s)| \le a_2|s|^{q_2} + b_2 \quad \forall s \in \mathbb{R}. \tag{1.43}$$

(a) *If $q_1q_2 = 1$ and $a_2 > 0$ is small then* (1.38) *admits a weak solution for any $\sigma > 0$, $\tilde\sigma > 0$, $\varrho > 0$, $\tilde\varrho > 0$.*
(b) *If $q_1q_2 < 1$ then* (1.38) *admits a weak solution for any $\sigma > 0$, $\tilde\sigma > 0$, $\varrho > 0$, $\tilde\varrho > 0$.*

III. Subcriticality and sublinearity. *Assume that $\Lambda_g < \infty$. In addition, assume that there exist $a_1 > 0$, $a_2 > 0$, $b_1 > 0$, $b_2 > 0$, $q_1 \in (1, p_\mu)$, $q_2 \in (0, 1]$, such that* (1.40) *and* (1.43) *hold.*
(a) *If $q_1q_2 > 1$ then* (1.38) *admits a weak solution for $b_1, b_2, \sigma, \tilde\sigma, \varrho, \tilde\varrho$ small enough.*
(b) *If $q_2p_\mu = 1$ and $a_2$ is mall enough then* (1.38) *admits a weak solution for any $\sigma > 0$, $\tilde\sigma > 0$, $\varrho > 0$, $\tilde\varrho > 0$.*
(c) *If $q_2p_\mu < 1$ then* (1.38) *admits a weak solution for every for any $\sigma > 0$, $\tilde\sigma > 0$, $\varrho > 0$, $\tilde\varrho > 0$.*

**Remark about elliptic equations and systems with weights.** We emphasize that Theorems B and C can be extended to the case of equations with weights of the form

$$-L_\mu u = \delta^\gamma u^p + \sigma\tau \quad \text{in } \Omega, \tag{1.44}$$

and Theorems E–H can be extended to the case of systems with weights of the form

$$\begin{cases} -L_\mu u = \epsilon\,\delta^\gamma g(v) + \sigma\tau & \text{in } \Omega, \\ -L_\mu v = \epsilon\,\delta^{\tilde\gamma}\tilde g(u) + \tilde\sigma\tilde\tau & \text{in } \Omega, \end{cases} \tag{1.45}$$

by using similar arguments. However, in order to avoid the complication of the proofs, we state and prove the results without weights.

The paper is organized as follows. In Section 2 we investigate properties of the boundary trace defined in Definition 1.2 and prove Propositions A and D. Theorems B and C are proved in Section 3 due to estimates on Green kernel, Martin kernel and the capacities $\mathrm{Cap}^{(p+1)\alpha}_{\mathbb{N}_{2\alpha,2},p'}$ and $\mathrm{Cap}^{\partial\Omega}_{1-\alpha+\frac{\alpha+1}{p},p'}$. In Section 4 sufficient conditions for the existence of weak solutions to elliptic systems with power source terms (1.32) (Theorems E and F) are obtained by combining the method in [6] and the capacity approach. Finally, in Section 5, we establish existence results for elliptic systems with more general nonlinearities (Theorems G and H) due to Schauder fixed point theorem.

**Notations.** Throughout this paper, $C$, $c$, $c'$,...denotes positive constants which may vary from one appearance to another. The notation $A \approx B$ means $c^{-1}B \leq A \leq cB$ for some constant $c > 1$ depending on some structural constant.

**Acknowledgment.** The first author has been supported by Fondecyt Grant 3140567 and by Millenium Nucleus CAPDE NC130017. The second author was supported by Fondecyt Grant 3160207.

## 2. Preliminaries

2.1. **Green kernel and Martin kernel.** Denote by $L^p_w(\Omega;\tau)$, $1 \leq p < \infty$, $\tau \in \mathfrak{M}^+(\Omega)$, the weak $L^p$ space (or Marcinkiewicz space) (see [17]). When $\tau = \delta^s dx$, for simplicity, we use the notation $L^p_w(\Omega;\delta^s)$. Notice that, for every $s > -1$,

$$L^p_w(\Omega;\delta^s) \subset L^r(\Omega;\delta^s) \quad \forall r \in [1,p). \tag{2.1}$$

Moreover for any $u \in L^p_w(\Omega;\delta^s)$ $(s > -1)$,

$$\int_{\{|u|\geq\lambda\}} \delta^s dx \leq \lambda^{-p} \|u\|^p_{L^p_w(\Omega;\delta^s)} \quad \forall \lambda > 0. \tag{2.2}$$

Let $G^\Omega_\mu$ and $K^\Omega_\mu$ be respectively the Green kernel and Martin kernel of $-L_\mu$ in $\Omega$ (see [14, 9]) for more details). We recall that

$$G^\Omega_\mu(x,y) \approx \min\left\{|x-y|^{2-N}, \delta(x)^\alpha\delta(y)^\alpha|x-y|^{2-N-2\alpha}\right\} \quad \forall x,y \in \Omega, x \neq y, \tag{2.3}$$

$$K^\Omega_\mu(x,y) \approx \delta(x)^\alpha|x-y|^{2-N-2\alpha} \quad \forall x \in \Omega,\, y \in \partial\Omega. \tag{2.4}$$

Finally, We denote by $\mathbb{G}_\mu$ and $\mathbb{K}_\mu$ be the corresponding Green operator and Martin operator (see [14, 9]), namely

$$\mathbb{G}_\mu[\tau](x) = \int_\Omega G_\mu(x,y)d\tau(y), \quad \forall \tau \in \mathfrak{M}(\Omega), \tag{2.5}$$

$$\mathbb{K}_\mu[\nu](x) = \int_{\partial\Omega} K_\mu(x,z)d\nu(z), \quad \forall \nu \in \mathfrak{M}(\partial\Omega). \tag{2.6}$$

Let us recall a result from [4] which will be useful in the sequel.

**Proposition 2.1.** ([4, Lemma 2.4]) *Let $\omega$ be a nonnegative bounded Radon measure in $D = \Omega$ or $\partial\Omega$ and $\eta \in C(\Omega)$ be a positive weight function. Let $H$ be a continuous nonnegative function on $\{(x,y) \in \Omega \times D :\ x \neq y\}$. For any $\lambda > 0$ we set*

$$A_\lambda(y) := \{x \in \Omega \setminus \{y\} :\ \ H(x,y) > \lambda\} \quad \textit{and} \quad m_\lambda(y) := \int_{A_\lambda(y)} \eta(x)dx.$$

*Suppose that there exist $C > 0$ and $k > 1$ such that $m_\lambda(y) \leq C\lambda^{-k}$ for every $\lambda > 0$. Then the operator*

$$\mathbb{H}[\omega](x) := \int_D H(x,y)d\omega(y)$$

*belongs to $L^k_w(\Omega;\eta)$ and*

$$||\mathbb{H}[\omega]||_{L^k_w(\Omega;\eta)} \leq (1 + \frac{Ck}{k-1})\omega(D).$$

By combining (2.3), (2.4) and the above Lemma we have the following result.

**Lemma 2.2.** *Let $\gamma \in \left(-\frac{\alpha N}{N+2\alpha-2}, \frac{\alpha N}{N-2}\right)$. Then there exists $C = C(N,\mu,\gamma,\Omega) > 0$ such that*

$$\sup_{\xi\in\Omega} \left|\left| \frac{G_\mu(\cdot,\xi)}{\delta(\xi)^\alpha} \right|\right|_{L_w^{\frac{N+\gamma}{N+\alpha-2}}(\Omega;\delta^\gamma)} < C. \tag{2.7}$$

*Proof.* Let $\xi \in \Omega$. We will apply Proposition 2.1 with $D = \Omega$, $\eta = \delta^\gamma$ with $\gamma > -1$, $\omega = \delta^\alpha\delta_\xi$, where $\delta_\xi$ is the Dirac measure concentrated at $\xi$, and

$$H(x,y) = \frac{G_\mu(x,y)}{\delta(y)^\alpha}.$$

Then

$$\mathbb{H}[\omega](x) = \int_\Omega \frac{G_\mu(x,y)}{\delta(y)^\alpha}\delta(y)^\alpha d\delta_\xi(y) = G_\mu(x,\xi).$$

From (2.3), there exists $C = C(N,\mu,\Omega)$ such that, for every $(x,y) \in \Omega\times\Omega$, $x \neq y$,

$$G_\mu(x,y) \leq C\delta(y)^\alpha|x-y|^{2-N-\alpha}, \tag{2.8}$$

$$G_\mu(x,y) \leq C\frac{\delta(y)^\alpha}{\delta(x)^\alpha}|x-y|^{2-N}, \tag{2.9}$$

$$G_\mu(x,y) \leq C\delta(x)^\alpha\delta(y)^\alpha|x-y|^{2-N-2\alpha}. \tag{2.10}$$

By (2.8), for any $x \in A_\lambda(y)$,

$$\lambda \leq C|x-y|^{2-N-\alpha}, \tag{2.11}$$

and form (2.9) and (2.10)

$$\delta(x)^\alpha \leq \frac{C}{\lambda}|x-y|^{2-N} \quad \text{and} \quad \delta(x)^\alpha \geq C\lambda|x-y|^{N+2\alpha-2} \tag{2.12}$$

We consider two cases: $\gamma \geq 0$ and $-1 < \gamma < 0$.
**Case 1:** $\gamma \geq 0$**.** Due to (2.11) and (2.12) we have

$$m_\lambda(y) = \int_{A_\lambda(y)} \delta(x)^\gamma dx \leq \int_{A_\lambda(y)} \left(\frac{C}{\lambda}|x-y|^{2-N}\right)^{\frac{\gamma}{\alpha}} dx \leq C\lambda^{-\frac{N+\gamma}{N+\alpha-2}},$$

with $\gamma < \frac{\alpha N}{N-2}$. Observe that $\omega(\Omega) = \delta(\xi)^\alpha$, by Proposition 2.1, we get

$$||G_\mu(\cdot,\xi)||_{L_w^{\frac{N+\gamma}{N+\alpha-2}}(\Omega;\,\delta^\gamma)} \leq C\delta(\xi)^\alpha.$$

This implies (2.7).

**Case 2:** $-1 < \gamma < 0$. By (2.11) and (2.12) we have

$$m_\lambda(y) = \int_{A_\lambda(y)} \delta(x)^\gamma dx \le \int_{A_\lambda(y)} (C\lambda|x-y|^{N+2\alpha-2})^{\frac{\gamma}{\alpha}} dx \le C\lambda^{-\frac{N+\gamma}{N+\alpha-2}},$$

with $\gamma > -\frac{\alpha N}{N+2\alpha-2}$. By arguing similarly as in Case 1, we get (2.7). □

**Lemma 2.3.** *Let $\gamma > -1$. Then there exists $C = C(N,\mu,\gamma,\Omega) > 0$ such that*

$$\sup_{\xi\in\partial\Omega} ||K_\mu(\cdot,\xi)||_{L_w^{\frac{N+\gamma}{N+\alpha-2}}(\Omega;\delta^\gamma)} < C.$$

*Proof.* Let $\xi \in \partial\Omega$. We will apply Proposition 2.1 with $D = \partial\Omega$, $\eta = \delta^\gamma$ with $\gamma > -1$ and $\omega = \delta_\xi$. The rest of the proof can be proceeded as in the proof of Lemma 2.2 and we omit it. □

In view of (2.1), Lemma 2.2 and Lemma 2.3, one can obtain easily the following proposition (see also [14, 18]).

**Proposition 2.4.** (i) *Let $\gamma \in (-\frac{\alpha N}{N+2\alpha-2}, \frac{\alpha N}{N-2})$. Then there exists a constant $c = c(N,\mu,\gamma,\Omega)$ such that*

$$\|\mathbb{G}_\mu[\tau]\|_{L_w^{\frac{N+\gamma}{N+\alpha-2}}(\Omega;\delta^\gamma)} \le c\,\|\tau\|_{\mathfrak{M}(\Omega;\delta^\alpha)} \quad \forall \tau \in \mathfrak{M}(\Omega;\delta^\alpha). \tag{2.13}$$

(ii) *Let $\gamma > -1$. Then there exists a constant $c = c(N,\mu,\gamma,\Omega)$ such that*

$$\|\mathbb{K}_\mu[\nu]\|_{L_w^{\frac{N+\gamma}{N+\alpha-2}}(\Omega;\delta^\gamma)} \le c\,\|\nu\|_{\mathfrak{M}(\partial\Omega)} \quad \forall \nu \in \mathfrak{M}(\partial\Omega). \tag{2.14}$$

2.2. **Boundary trace.** In this section we study properties of the boundary trace in connection with of $L_\mu$ harmonic functions. In particular, we show that, when $\mu < C_H(\Omega)$, the boundary trace defined in Definition 1.2 coincides the notion of normalized boundary trace introduced in Definition 1.1). To this end, we will examine that $\operatorname{tr}(\mathbb{G}_\mu[\tau]) = 0$ for every $\tau \in \mathfrak{M}(\Omega;\delta^\alpha)$ and $\operatorname{tr}(\mathbb{K}_\mu[\nu]) = \nu$ for every $\nu \in \mathfrak{M}(\partial\Omega)$. These results are proved below, based on a combination of the ideas in [9] and [14]. It is worth emphasizing that the below results are valid for $\mu \in (0, \frac{1}{4}]$ (under the condition that the first eigenvalue $\lambda_\mu$ of $-L_\mu$ is positive).

**Proposition 2.5.** *Let $\tau \in \mathfrak{M}(\Omega;\delta^\alpha)$ and $u = \mathbb{G}_\mu[\tau]$. Then $\operatorname{tr}(u) = 0$.*

*Proof.* First we assume that $\tau$ is nonnegative. Let $\{\Omega_n\}$ be a smooth exhaustion of $\Omega$ and for each $n$, let $\omega_{\Omega_n}^{x_0}$ be the $L_\mu^{\Omega_n}$ harmonic measure on $\partial\Omega_n$. Then $u$ satisfies

$$\begin{cases} -L_\mu u = \tau & \text{in } \Omega_n \\ \quad u = u & \text{on } \partial\Omega_n. \end{cases} \tag{2.15}$$

Thus

$$u = \mathbb{G}_\mu^{\Omega_n}[\tau] + \mathbb{K}_\mu^{\Omega_n}[u] = \mathbb{G}_\mu^{\Omega_n}[\tau] + \int_{\partial\Omega_n} u d\omega_{\Omega_n}^{x_0}. \tag{2.16}$$

This, joint with $\mathbb{G}_\mu^{\Omega_n}[\tau] \uparrow \mathbb{G}_\mu[\tau]$ as $n \to \infty$, ensures

$$\lim_{n\to\infty} \int_{\partial\Omega_n} u d\omega_{\Omega_n}^{x_0} = 0,$$

namely $\operatorname{tr}(u) = 0$.

In the general case, the result follows from the linearity property of the problem. □

The next result shows that the boundary trace of $L_\mu$ harmonic function can be achieved in a *dynamic way.*

**Proposition 2.6.** [9, Proposition 2.34] *Let $x_0 \in \Omega_1$ and $\mu \in \mathfrak{M}(\partial\Omega)$. Put*

$$v(x) := \int_{\partial\Omega} K_\mu(x,y)d\nu(y),$$

*then for every $\zeta \in C(\overline{\Omega})$,*

$$\lim_{n\to\infty}\int_{\partial\Omega_n} \zeta v d\omega_{\Omega_n}^{x_0} = \int_{\partial\Omega}\zeta d\nu. \tag{2.17}$$

Also we have the following representation formula for $L_\mu$ harmonic functions.

**Proposition 2.7.** [9, Theorem 2.33] *Let $u$ be a positive $L_\mu$ harmonic in $\Omega$. Then $u \in L^1(\Omega;\delta^\alpha)$ and there exists a unique Radon measure $\nu$ on $\partial\Omega$ such that*

$$u(x) = \int_{\partial\Omega} K_\mu(x,y)d\nu(y). \tag{2.18}$$

In the following proposition, we study the boundary trace of $L_\mu$ subharmonic functions.

**Proposition 2.8.** *Let $w$ be a nonnegative $L_\mu$ subharmonic function. If $w$ is dominated by an $L_\mu$ superharmonic function then $L_\mu w \in \mathfrak{M}^+(\Omega;\delta^\alpha)$ and $w$ has a boundary trace $\nu \in \mathfrak{M}(\partial\Omega)$. In addition, if* $\mathrm{tr}\,(w) = 0$ *then $w = 0$.*

*Proof.* By proceeding as in the proof of [14, Proposition 2.14] and using Proposition 2.7, we obtain the desired result. □

**Proposition 2.9.** *Let $w$ be a nonnegative $L_\mu$ subharmonic function. If $w$ has a boundary trace then it is dominated by an $L_\mu$ harmonic function.*

*Proof.* The proof is similar to that of Proposition 2.20 in [14]. For the sake of convenience we give it below. Let $\{\Omega_n\}$ be as in the proof of Proposition 2.5 and fix $x_0 \in \Omega_1$. For any $x \in \Omega$, set

$$u_n(x) = \int_{\partial\Omega_n} w d\omega_{\Omega_n}^x,$$

then $u_n$ is $L_\mu^{\Omega_n}$ harmonic function with boundary trace $w$. Furthermore, by the maximum principle we have that $w \le u_n$ in $\Omega_n$. Let $\nu \in \mathfrak{M}(\partial\Omega)$ be such that

$$\lim_{n\to\infty}\int_{\partial\Omega_n}\zeta w d\omega_{\Omega_n}^{x_0} = \int_{\partial\Omega}\zeta d\nu \quad \forall \zeta \in C(\overline{\Omega}). \tag{2.19}$$

Then

$$u_n(x_0) = \int_{\partial\Omega_n} w d\omega_{\Omega_n}^{x_0} \to \int_{\partial\Omega} d\nu.$$

We infer from Harnack inequality that $\{u_n\}$ is locally uniformly bounded and hence there exists an $L_\mu$ harmonic function $u$ such that $u_n \to u$ locally uniformly in $\Omega$. By Proposition 2.8, there exists a nonnegative measure $\tau \in \mathfrak{M}^+(\Omega;\delta^\alpha)$ such that

$$w = -\mathbb{G}_\mu[\tau] + \mathbb{K}_\mu[\nu].$$

On the other hand,

$$w = -\mathbb{G}_\mu^{\Omega_n}[\tau] + u_n \to -\mathbb{G}_\mu[\tau] + u,$$

locally uniformly in $\Omega$. Thus we can deduce that $u = \mathbb{K}_\mu[\nu]$ and the result follows. □

**Proposition 2.10.** *Let $u$ be a nonnegative $L_\mu$ superharmonic function. Then there exist $\nu \in \mathfrak{M}^+(\partial\Omega)$ and $\tau \in \mathfrak{M}^+(\Omega;\delta^\alpha)$ such that*

$$u = \mathbb{G}_\mu[\tau] + \mathbb{K}_\mu[\nu].$$

*Proof.* Let $\Omega_n$ and $\omega_{\Omega_n}^{x_0}$ be as in the proof of Proposition 2.5. Since $u$ is $L_\mu$ superharmonic function there exists a nonnegative Radon measure in $\Omega$ such that

$$-L_\mu u = \tau \quad \text{in } \Omega$$

in the sense of distributions. Note that $u$ is the unique solution of

$$\left\{ \begin{aligned} -L_\mu w &= \tau && \text{in } \ \Omega_n \\ w &= u && \text{on } \ \partial\Omega_n. \end{aligned} \right. \tag{2.20}$$

Therefore

$$u = \mathbb{G}_\mu^{\Omega_n}[\tau] + \mathbb{K}_\mu^{\Omega_n}[u]. \tag{2.21}$$

Set $w_n = \mathbb{K}_\mu^{\Omega_n}[u]$. Since $\tau \geq 0$, by the above quality, we have $0 \leq w_n(x) \leq u(x)$. Thus by the Harnack inequality, $w_n \to w$ locally uniformly in $\Omega$. Furthermore, $w$ is an $L_\mu$ harmonic function in $\Omega$ and by Proposition 2.18 there exists $\nu \in \mathfrak{M}^+(\partial\Omega)$ such that

$$w = \mathbb{K}_\mu[\nu]. \tag{2.22}$$

Now since $G_\mu^{\Omega_n} \uparrow G_\mu$ as $n \to \infty$, we deduce from (2.21) and (2.22) that

$$u = \mathbb{G}_\mu^{\Omega_n}[\tau] + \mathbb{K}_\mu^{\Omega_n}[u] \to \mathbb{G}_\mu[\tau] + \mathbb{K}_\mu[\nu].$$

Since

$$G_\mu(x,y) \geq c(x,\mu,N)\delta(y)^\alpha,$$

we can easily prove that $\tau \in \mathfrak{M}^+(\Omega;\delta^\alpha)$ which completes the proof. ☐

The above results enable to study the boundary value problem for the linear equation

$$\left\{ \begin{aligned} -L_\mu u &= \tau && \text{in } \ \Omega, \\ \operatorname{tr}(u) &= \nu. \end{aligned} \right. \tag{2.23}$$

**Definition 2.11.** *Let $(\tau,\nu) \in \mathfrak{M}(\Omega;\delta^\alpha) \times \mathfrak{M}(\partial\Omega)$. We say that $u$ is a weak solution of (2.23) if $u \in L^1(\Omega;\delta^\alpha)$ and*

$$-\int_\Omega u L_\mu \zeta \, dx = \int_\Omega \zeta \, d\tau - \int_\Omega \mathbb{K}_\mu[\nu] L_\mu \zeta \, dx \qquad \forall \zeta \in \mathbf{X}_\mu(\Omega), \tag{2.24}$$

**Proposition 2.12.** *For any $(\tau,\nu) \in \mathfrak{M}(\Omega;\delta^\alpha) \times \mathfrak{M}(\partial\Omega)$ there exists a unique weak solution of (2.23). Moreover*

$$u = \mathbb{G}_\mu[\tau] + \mathbb{K}_\mu[\nu], \tag{2.25}$$

$$\|u\|_{L^1(\Omega;\delta^\alpha)} \leq c(\|\tau\|_{\mathfrak{M}(\Omega;\delta^\alpha)} + \|\nu\|_{\mathfrak{M}(\partial\Omega)}). \tag{2.26}$$

*In addition, for any $\zeta \in \mathbf{X}_\mu(\Omega)$, $\zeta \geq 0$,*

$$-\int_\Omega |u| L_\mu \zeta \, dx \leq \int_\Omega \zeta \operatorname{sign}(u) \, d\tau - \int_\Omega \mathbb{K}_\mu[|\nu|] L_\mu \zeta \, dx, \tag{2.27}$$

*and*

$$-\int_\Omega u_+ L_\mu \zeta \, dx \leq \int_\Omega \zeta \operatorname{sign}_+(u) \, d\tau - \int_\Omega \mathbb{K}_\mu[\nu_+] L_\mu \zeta \, dx. \tag{2.28}$$

*Proof.* The proof is similar to that of [9, Proposition 3.2] and we omit it. □

**Remark 2.1.** *If $h \in L^1(\partial\Omega, d\omega_\Omega^{x_0})$ is the boundary value of* (2.23)*, the above Proposition is valid for $d\nu = hd\omega_\Omega^{x_0}$.*

**Proposition 2.13.** *(i) For $\tau \in \mathfrak{M}(\Omega;\delta^\alpha)$, $\mathrm{tr}\,(\mathbb{G}_\mu[\tau]) = 0$ and for $\nu \in \mathfrak{M}(\partial\Omega)$, $\mathrm{tr}\,(\mathbb{K}_\mu[\nu]) = \nu$.*

*(ii) Let $w$ be a nonnegative $L_\mu$ subharmonic function in $\Omega$. Then $w$ is dominated by an $L_\mu$ superharmonic function if and only if $w$ has a boundary trace $\nu \in \mathfrak{M}(\partial\Omega)$. Moreover, if $w$ has a boundary trace then $L_\mu w \in \mathfrak{M}^+(\Omega;\delta^\alpha)$. If, in addition, if $\mathrm{tr}\,(w) = 0$ then $w = 0$.*

*(iii) Let $u$ be a nonnegative $L_\mu$ superharmonic function. Then there exist $\nu \in \mathfrak{M}^+(\partial\Omega)$ and $\tau \in \mathfrak{M}^+(\Omega,\delta^\alpha)$ such that* (2.25) *holds.*

*(iv) Let $(\tau,\nu) \in \mathfrak{M}(\Omega;\delta^\alpha) \times \mathfrak{M}(\partial\Omega)$. Then there exists a unique weak solution $u$ of* (2.23)*. The solution is given by* (2.25)*. Moreover, there exists $c = c(N,\mu,\Omega)$ such that* (2.26) *holds.*

*Proof.* Statement (i) follows from Proposition 2.5 and Proposition 2.6. Statement (ii) can be deduced from Proposition 2.8 and Proposition 2.9. Statement (iii) follows from Proposition 2.10. Finally statement (iv) is obtained due to Proposition 2.12. □

**Proof of Proposition A.** We infer from [9] that $(i) \Longleftrightarrow (ii)$. By an argument similar to that of the proof of [18, Theorem B], we deduce that $(ii) \Longleftrightarrow (iii)$. □

For $\beta > 0$, put

$$\Omega_\beta := \{x \in \Omega : \delta(x) < \beta\},\ D_\beta := \{x \in \Omega : \delta(x) > \beta\},\ \Sigma_\beta := \{x \in \Omega : \delta(x) = \beta\}. \tag{2.29}$$

**Lemma 2.14.** *There exists $\beta_* > 0$ such that for every point $x \in \overline{\Omega}_{\beta_*}$, there exists a unique point $\sigma_x \in \partial\Omega$ such that $x = \sigma_x - \delta(x)\mathbf{n}_{\sigma_x}$. The mappings $x \mapsto \delta(x)$ and $x \mapsto \sigma_x$ belong to $C^2(\overline{\Omega}_{\beta_*})$ and $C^1(\overline{\Omega}_{\beta_*})$ respectively. Moreover, $\lim_{x\to\sigma(x)} \nabla\delta(x) = -\mathbf{n}_{\sigma_x}$.*

**Proof of Proposition D.**
**(iii) $\Longrightarrow$ (ii).** Assume $(u,v)$ is a distribution solution of (1.30). Put $\omega := g(v)$ and denote $\omega_\beta := \omega|_{D_\beta}$, $\tau_\beta := \tau|_{D_\beta}$ and $\lambda_\beta := u|_{\Sigma_\beta}$ for $\beta \in (0,\beta_*)$. Consider the boundary value problem

$$-L_\mu w = \omega_\beta + \tau_\beta \quad \text{in } D_\beta, \qquad w = \lambda_\beta \quad \text{on } \Sigma_\beta.$$

This problem admits a unique solution $w_\beta$ (see [9]). Therefore $w_\beta = u|_{D_\beta}$. We have

$$u|_{D_\beta} = w_\beta = \mathbb{G}_\mu^{D_\beta}[\omega_\beta] + \mathbb{G}_\mu^{D_\beta}[\tau_\beta] + \mathbb{P}_\mu^{D_\beta}[\lambda_\beta]$$

where $\mathbb{G}_\mu^{D_\beta}$ and $\mathbb{P}_\mu^{D_\mu}$ are respectively Green kernel and Poisson kernel of $-L_\mu$ in $D_\beta$.

It follows that

$$\left|\int_{D_\beta} G_\mu^{D_\beta}(\cdot,y)g(v(y))\,dy\right| = \left|\mathbb{G}_\mu^{D_\beta}[\tau_\beta]\right| \leq \left|u|_{D_\beta}\right| + \left|\mathbb{G}_\mu^{\Omega}[\tau]\right| + \left|\mathbb{K}_\mu^{\Omega}[\nu]\right|.$$

Letting $\beta \to 0$, we get

$$\left|\int_\Omega G_\mu(\cdot,y)g(v(y))\,dy\right| < \infty. \tag{2.30}$$

Fix a point $x_0 \in \Omega$. Keeping in mind that $G_\mu(x_0,y) \approx \delta(y)^\alpha$ for every $y \in \Omega_{\beta_*}$, we deduce from (2.30) that $g(v) \in L^1(\Omega;\delta^\alpha)$. Similarly, one can show that $\tilde{g}(u) \in L^1(\Omega;\delta^\alpha)$. Thanks to Proposition 2.13 (v), we obtain (1.29).

**(ii) $\Longrightarrow$ (iii).** Assume $u$ and $v$ are functions such that $\tilde{g}(u) \in L^1(\Omega;\delta^\alpha)$, $g(v) \in L^1(\Omega;\delta^\alpha)$ and (1.29) holds. By Proposition 2.13 (i) $L_\mu\mathbb{K}_\mu[\nu] = L_\mu\mathbb{K}_\mu[\tilde{\nu}] = 0$, which implies that $(u,v)$ is a solution of (1.30). On the other hand, since $\tilde{g}(u) \in L^1(\Omega;\delta^\alpha)$ and $g(v) \in L^1(\Omega;\delta^\alpha)$, we deduce from Proposition 2.13 (ii) that $\mathrm{tr}\,(\mathbb{G}_\mu[\tilde{g}(u)]) = \mathrm{tr}\,(\mathbb{G}_\mu[g(v)]) = 0$. Consequently, $\mathrm{tr}\,(u) = \mathrm{tr}\,(\mathbb{K}_\mu[\nu]) = \nu$ and $\mathrm{tr}\,(v) = \mathrm{tr}\,(\mathbb{K}_\mu[\tilde{\nu}]) = \tilde{\nu}$.

**(iii) $\Longrightarrow$ (i).** Assume $(u,v)$ is a positive solution of (1.30) in the sense of distributions. From the implication (iii) $\Longrightarrow$ (ii), we deduce that $u \in L^1(\Omega;\delta^\alpha)$, $v \in L^1(\Omega;\delta^\alpha)$, $\tilde{g}(u) \in L^1(\Omega;\delta^\alpha)$ and $g(v) \in L^1(\Omega;\delta^\alpha)$. Hence, by Proposition 2.13, (1.28) holds for every $\phi \in \mathbf{X}_\mu(\Omega)$.

**(i) $\Longrightarrow$ (iii).** This implication follows straight foward from Proposition 2.13. □

## 3. The scalar problem

3.1. **Concavity properties and Green properties.** Here we give some concavity lemmas that will be employed in the sequel.

**Proposition 3.1.** *Let $\varphi \in L^1(\Omega;\delta^\alpha)$, $\varphi \geq 0$ and $\tau \in \mathfrak{M}^+(\Omega;\delta^\alpha)$. Set*

$$w := \mathbb{G}_\mu[\varphi + \tau] \quad \text{and} \quad \psi = \mathbb{G}_\mu[\tau].$$

*Let $\phi$ be a concave nondecreasing $C^2$ function on $[0,\infty)$, such that $\phi(1) \geq 0$. Then $\phi'(w/\psi)\varphi \in L^1(\Omega;\delta^\alpha)$ and the following holds in the weak sense in $\Omega$*

$$-L_\mu(\psi\phi(w/\psi)) \geq \phi'(w/\psi)\varphi.$$

*Proof.* Let $\{\varphi_n\}, \tau_n \in C^\infty(\overline{\Omega})$ such that $\varphi_n \to \varphi$ in $L^1(\Omega,\delta^\alpha)$ and $\tau_n \rightharpoonup \tau$. Set $w_n := \mathbb{G}_\mu[\varphi_n + \tau_n]$ and $\psi_n = \mathbb{G}_\mu[\tau_n]$. Since $w_n \geq \psi_n > 0$ for any $n \geq n_0$ for some $n_0 \in \mathbb{N}$, we have by straightforward calculations

$$-\Delta\left(\psi_n\phi(\frac{w_n}{\psi_n})\right) = (-\Delta\psi_n)\left(\phi(\frac{w_n}{\psi_n}) - \frac{w_n}{\psi_n}\phi'(\frac{w_n}{\psi_n})\right) + (-\Delta w_n)\phi'(\frac{w_n}{\psi_n}) - \psi_n\phi''(\frac{w_n}{\psi_n})\left|\nabla\left(\frac{w_n}{\psi_n}\right)\right|^2.$$

Now note that, since $\phi' \geq 0$, we have

$$(-\Delta w_n)\phi'(\frac{w_n}{\psi_n}) \geq \phi'(\frac{w_n}{\psi_n})\left(-\Delta\psi_n - \mu\frac{\psi_n}{\delta^2} + \mu\frac{w_n}{\delta^2} + \varphi_n\right).$$

This, together with the fact that $\phi(t) - t\phi'(t) + \phi'(t) \geq 0$ for any $t \geq 1$, implies

$$\begin{aligned}
&(-\Delta\psi_n)\left(\phi(\frac{w_n}{\psi_n}) - \frac{w_n}{\psi_n}\phi'(\frac{w_n}{\psi_n})\right) + (-\Delta w_n)\phi'(\frac{w_n}{\psi_n}) \\
&\geq (-\Delta\psi_n)\left(\phi(\frac{w_n}{\psi_n}) - \frac{w_n}{\psi_n}\phi'(\frac{w_n}{\psi_n}) + \phi'(\frac{w_n}{\psi_n})\right) + \phi'(\frac{w_n}{\psi_n})\left(-\mu\frac{\psi_n}{\delta^2} + \mu\frac{w_n}{\delta^2} + \varphi_n\right) \\
&\geq \mu\frac{\psi_n}{\delta^2}\left(\phi(\frac{w_n}{\psi_n}) - \frac{w_n}{\psi_n}\phi'(\frac{w_n}{\psi_n}) + \phi'(\frac{w_n}{\psi_n})\right) + \phi'(\frac{w_n}{\psi_n})\left(-\mu\frac{\psi_n}{\delta^2} + \mu\frac{w_n}{\delta^2} + \varphi_n\right) \\
&= \frac{\mu}{\delta^2}\psi_n\phi(\frac{w_n}{\psi_n}) + \phi'(\frac{w_n}{\psi_n})\varphi_n.
\end{aligned}$$

Thus we have proved

$$-L_\mu\left(\psi_n\phi(\frac{w_n}{\psi_n})\right) \geq \phi'(\frac{w_n}{\psi_n})\varphi_n.$$

Also

$$\psi_n\phi(\frac{w_n}{\psi_n}) \leq \psi_n(\phi(0) + \phi'(0)\frac{w_n}{\psi_n}) \leq C(\psi_n + w_n)$$

and

$$-\int_\Omega \psi_n\phi(\frac{w_n}{\psi_n})L_\mu\xi\,dx \geq \int_\Omega \phi'(\frac{w_n}{\psi_n})\varphi_n\xi\,dx \quad \forall\xi\in\mathbf{X}_\mu(\Omega).$$

By passing to the limit with Lebesgue theorem and Fatou lemma, we complete the proof. ☐

In the next Lemma we will prove the 3-$G$ inequality which will be useful later.

**Lemma 3.2.** *There exists a positive constant $C=C(N,\mu,\Omega)$ such that*

$$\frac{G_\mu(x,y)G_\mu(y,z)}{G_\mu(x,z)} \leq C\left(\frac{\delta(y)^\alpha}{\delta(x)^\alpha}G_\mu(x,y)+\frac{\delta(y)^\alpha}{\delta(z)^\alpha}G_\mu(y,z)\right) \quad \forall(x,y,z)\in\Omega\times\Omega\times\Omega. \tag{3.1}$$

*Proof.* It follows from (2.3) and the inequality $|\delta(x)-\delta(y)|\leq|x-y|$ that

$$\begin{aligned} G_\mu(x,y) &\approx \min\left\{|x-y|^{2-N},\delta(x)^\alpha\delta(y)^\alpha|x-y|^{2-2\alpha-N}\right\} \\ &\approx |x-y|^{2-N}\delta(x)^\alpha\delta(y)^\alpha\left(\max\left\{\delta(x)^\alpha\delta(y)^\alpha,|x-y|^{2\alpha}\right\}\right)^{-1} \\ &\approx |x-y|^{2-N}\delta(x)^\alpha\delta(y)^\alpha\left(\max\left\{\delta(x),\delta(y),|x-y|\right\}\right)^{-2\alpha} \\ &= \delta(x)^\alpha\delta(y)^\alpha N_{2\alpha,2}(x,y), \quad \forall x,y\in\Omega,\ x\neq y, \end{aligned}$$

where $N_{2\alpha,2}(x,y)$ is defined in (1.23) with $a=2\alpha$ and $\beta=2$. By [5, Lemma 2.2] we deduce that there exists a positive constant $C=C(N,\mu,\Omega)$ such that

$$\frac{1}{N_{2\alpha,2}(x,z)} \leq C\left(\frac{1}{N_{2\alpha,2}(x,y)}+\frac{1}{N_{2\alpha,2}(y,z)}\right). \tag{3.2}$$

From (3.2) we can easily obtain (3.1). ☐

**Lemma 3.3.** *Let $0<p<p_\mu$ and $\tau\in\mathfrak{M}^+(\Omega;\delta^\alpha)$. Then there is a constant $C=C(N,\mu,p,\tau,\Omega)>0$ such that (1.19) holds.*

*Proof.* First we assume that $p>1$. By (2.13) we have that $\mathbb{G}_\mu[\tau]^p\in L^1(\Omega;\delta^\alpha)$. We write

$$\mathbb{G}_\mu[\tau](y)=\int_\Omega G_\mu(y,z)d\tau(z)=\int_\Omega\frac{G_\mu(y,z)}{\delta(z)^\alpha}\delta(z)^\alpha d\tau(z),$$

thus

$$\mathbb{G}_\mu[\tau](y)^p\leq C\int_\Omega\delta(z)^\alpha\left(\frac{G_\mu(y,z)}{\delta(z)^\alpha}\right)^p d\tau(z).$$

Consequently,

$$\mathbb{G}_\mu[\mathbb{G}_\mu[\tau]^p](x)\leq C\int_\Omega\int_\Omega G_\mu(x,y)G_\mu(y,z)^p\delta(z)^{\alpha(1-p)}d\tau(z)dy. \tag{3.3}$$

Also by (3.1) we obtain

$$\begin{aligned} &\int_\Omega\int_\Omega G_\mu(x,y)G_\mu(y,z)^p\delta(z)^{\alpha(1-p)}d\tau(z)dy \\ &\leq C\int_\Omega G_\mu(x,z)\int_\Omega\delta(y)^\alpha\left(\left(\frac{G_\mu(x,y)}{\delta(x)^\alpha}\right)\left(\frac{G_\mu(y,z)}{\delta(z)^\alpha}\right)^{p-1}+\left(\frac{G_\mu(y,z)}{\delta(z)^\alpha}\right)^p\right)dyd\tau(z) \\ &\leq C\int_\Omega G_\mu(x,z)\int_\Omega\delta(y)^\alpha\left(\left(\frac{G_\mu(x,y)}{\delta(x)^\alpha}\right)^p+\left(\frac{G_\mu(y,z)}{\delta(z)^\alpha}\right)^p\right)dyd\tau(z), \qquad (3.4) \end{aligned}$$

where in the last inequality we have used the Hölder inequality. By (3.3), (3.4) and Lemma 2.2 we derive that

$$\mathbb{G}_\mu[\mathbb{G}_\mu[\tau]^p](x) \le C \int_\Omega G_\mu(x,z)d\tau(z).$$

Note that the above argument is still valid for $p=1$.

If $0 \le p < 1$ then

$$\mathbb{G}_\mu[\mathbb{G}_\mu[\tau]^p] \le C(\mathbb{G}_\mu[1] + \mathbb{G}_\mu[\mathbb{G}_\mu[\tau]]).$$

By combining the case $p=1$ and the estimate $\mathbb{G}_\mu[1] \le C\mathbb{G}_\mu[\tau]$, we obtain (1.19). □

Actually (1.19) is a sufficient condition for the existence of weak solution of

$$\begin{cases} -L_\mu u = u^p + \sigma\tau & \text{in } \Omega, \\ \mathrm{tr}\,(u) = 0 & \text{on } \partial\Omega. \end{cases} \tag{3.5}$$

**Proposition 3.4.** *Let $0 < p \neq 1$, $\sigma > 0$ and $\tau \in \mathfrak{M}^+(\Omega;\delta^\alpha)$. Assume that there exists a positive constant $C$ such that (1.19) holds. Then problem (3.5) admits a weak solution $u$ satisfying*

$$\mathbb{G}_\mu[\sigma\tau] \le u \le C\mathbb{G}_\mu[\sigma\tau] \quad \text{a.e. in } \Omega, \tag{3.6}$$

*with another constant $C>0$, for any $\sigma > 0$ small enough if $p>1$, for any $\sigma > 0$ if $p < 1$.*

*Proof.* We adapt the idea in the proof of [6, Theorem 3.4]. Put $w := A\mathbb{G}_\mu[\sigma\tau]$ where $A > 0$ will be determined later. By (1.19),

$$\mathbb{G}_\mu[w^p + \sigma\tau] \le (CA^p\sigma^{p-1} + 1)\mathbb{G}_\mu[\sigma\tau] \quad \text{in } \Omega.$$

Therefore we deduce that $w \ge \mathbb{G}_\mu[w^p + \sigma\tau]$ as long as

$$CA^p\sigma^{p-1} + 1 \le A. \tag{3.7}$$

If $p > 1$ then (3.7) holds if we choose $A > 1$ and then choose $\sigma > 0$ small enough. If $p \in (0,1)$ then (3.7) holds if we choose $\sigma > 0$ arbitrary and then choose $A > 0$ large enough.

Next put $u_0 := \mathbb{G}_\mu[\sigma\tau]$ and $u_{n+1} := \mathbb{G}_\mu[u_n^p + \sigma\tau]$. It is clear that $\{u_n\}$ is increasing and $u_n \le w$ in $\Omega$ for all $n$. Since (1.19) holds, $w^p \in L^1(\Omega;\delta^\alpha)$. Consequently, by monotone convergence theorem, there exists a function $u \in L^p(\Omega;\delta^\alpha)$ such that $u_n^p \to u^p$ in $L^1(\Omega;\delta^\alpha)$. It is easy to see that $u$ is a solution of (3.5) satisfying (3.6). □

Estimate (1.19) is also a necessary condition for the existence of weak solution of (3.5).

**Proposition 3.5.** *Let $p > 1$, $\sigma > 0$ and $\tau \in \mathfrak{M}^+(\Omega;\delta^\alpha)$. Assume that problem (3.5) admits a weak solution. Then (1.19) holds with $C = \frac{1}{p-1}$.*

*Proof.* We adapt the argument used in the proof of [6, Proposition 3.5]. Assume (3.5) has a solution $u \in L^p(\Omega;\delta^\alpha)$ and assume $\sigma = 1$. By applying Proposition 3.1 with $\varphi$ replaced by $u^p$ and with

$$\phi(s) = \begin{cases} (1 - s^{1-p})/(p-1) & \text{if } s \ge 1, \\ s - 1 & \text{if } s < 1, \end{cases}$$

we get (1.19) with $C = \frac{1}{p-1}$. □

**Proposition 3.6.** *Let $0 < p < p_\mu$, $\sigma > 0$ and $\tau \in \mathfrak{M}^+(\Omega;\delta^\alpha)$. Then there exists a positive constant $C = C(N,\mu,\Omega,\sigma,\tau)$ such that for any weak solution $u$ of (3.5) there holds*

$$\mathbb{G}_\mu[\sigma\tau] \le u \le C(\mathbb{G}_\mu[\sigma\tau] + \delta^\alpha) \quad \text{a.e. in } \Omega. \tag{3.8}$$

*Proof.* We follows the idea in the proof of [6, Theorem 3.6]. We may assume that $\sigma = 1$. If $0 \le p < 1$, then

$$u = \mathbb{G}_\mu[u^p + \tau] \le C(\mathbb{G}_\mu[1] + \mathbb{G}_\mu[u] + \mathbb{G}_\mu[\tau]).$$

Since $\mathbb{G}_\mu[1] \le C\delta^\alpha$ a.e. in $\Omega$, we obtain

$$u \le C(\delta^\alpha + \mathbb{G}_\mu[u] + \mathbb{G}_\mu[\tau]) \quad \text{a.e. in } \Omega.$$

Therefore it is sufficient to deal with the term $\mathbb{G}_\mu[u]$ and we may assume that $p \ge 1$. Set

$$u_1 := u - \mathbb{G}_\mu[\tau] = \mathbb{G}_\mu[u^p],$$

hence $u = u_1 + \mathbb{G}_\mu[\tau]$. Since $u \in L^p(\Omega; \delta^\alpha)$ (by assumption), it follows that $u + \tau \in \mathfrak{M}^+(\Omega; \delta^\alpha)$, therefore by (2.13), $u \in L^s(\Omega; \delta^\alpha)$, for all $1 \le s < p_\mu$. Thus there exists $k_0 > 1$ such that $u^p \in L^{k_0}(\Omega; \delta^\alpha)$.

Let $p < s < p_\mu$. By Hölder inequality we obtain

$$\begin{aligned} u_1(x)^{k_0 s} = \mathbb{G}_\mu[u^p](x)^{k_0 s} &= \left( \int_\Omega \frac{G_\mu(x,y)}{\delta(y)^\alpha} \delta(y)^\alpha u(y)^p dy \right)^{k_0 s} \\ &\le \left( \int_\Omega \frac{G_\mu(x,y)}{\delta(y)^\alpha} \delta(y)^\alpha u(y)^{k_0 p} dy \right)^s \left( \int_\Omega \frac{G_\mu(x,y)}{\delta(y)^\alpha} \delta(y)^\alpha dy \right)^{\frac{(k_0-1)s}{k_0}} \\ &\le C \int_\Omega \left( \frac{G_\mu(x,y)}{\delta(y)^\alpha} \right)^s \delta(y)^\alpha u(y)^{k_0 p} dy. \end{aligned}$$

This, joint with Lemma 2.2, yields

$$\int_\Omega u_1(y)^{k_0 s} \delta(y)^\alpha dy \le C \int_\Omega u(y)^{k_0 p} \delta(y)^\alpha \int_\Omega \left( \frac{G_\mu(x,y)}{\delta(y)^\alpha} \right)^s \delta(x)^\alpha dx dy < c.$$

Since $u^p \le C(u_1^p + \mathbb{G}_\mu[\tau]^p)$, by Lemma 3.3 we have

$$u \le C\left(\mathbb{G}_\mu[\mathbb{G}_\mu[\tau]^p] + u_2\right) + \mathbb{G}_\mu[\tau] \le C(\mathbb{G}_\mu[\tau] + u_2),$$

where $u_2 := \mathbb{G}_\mu[u_1^p]$. Note that $u_2 \in L^{\frac{k_0 s^2}{p^2}}(\Omega; \delta^\alpha)$.

By induction we define $u_n := \mathbb{G}_\mu[u_{n-1}^p]$ and we have $u \le C(\mathbb{G}_\mu[\tau] + u_n)$, $u_n^p \in L^{s_n}(\Omega; \delta^\alpha)$ with $s_n = \frac{k_0 s^n}{p^n}$. Since $s_n \to \infty$, by [17, Lemma 2.3.2] we have for $1 < s < p_\mu$,

$$\begin{aligned} u_n &\le C \int_\Omega |x-y|^{2-\alpha-N} u_{n-1}^p \delta(y)^\alpha dy \\ &\le C \left( \int_\Omega |x-y|^{(2-\alpha-N)s} \delta(y)^\alpha dy + \int_\Omega |x-y|^{2-\alpha-N} u_{n-1}^{\frac{ps}{s-1}} \delta(y)^\alpha dy \right) \\ &\le C', \end{aligned}$$

for $n$ large enough. Therefore we obtain $u \le C(\mathbb{G}_\mu[\tau] + 1)$, which implies $u \le C(\mathbb{G}_\mu[\tau] + \mathbb{G}_\mu[1])$ with another $C > 0$. This, together with the inequality $\mathbb{G}_\mu[1] \le C\delta^\alpha$, implies (3.8). ☐

### 3.2. New Green properties.

**Lemma 3.7.** *Let $0 < p < p_\mu$, $\tau \in \mathfrak{M}^+(\Omega; \delta^\alpha)$. Let $s$ be such that*

$$\max(0, p - p_\mu + 1) < s \le 1. \tag{3.9}$$

*Then there exists a constant $C > 0$ such that*

$$\mathbb{G}_\mu[\mathbb{G}_\mu[\tau]^p] \le C \mathbb{G}_\mu[\tau]^s \quad \text{a.e. in } \Omega. \tag{3.10}$$

*Proof.* First we assume that $p > 1$. In view of the proof of Lemma 3.3, we have

$$\mathbb{G}_\mu[\mathbb{G}_\mu[\tau]^p](x) \le C\int_\Omega\int_\Omega G_\mu(x,y)G_\mu(y,z)^p\delta(z)^{\alpha(1-p)}d\tau(z)dy$$

$$= C\int_\Omega\int_\Omega G_\mu(x,y)^{1-s}G_\mu(x,y)^sG_\mu(y,z)^s\left(\frac{G_\mu(y,z)}{\delta^\alpha(z)}\right)^{p-s}\delta(z)^{\alpha(1-s)}d\tau(z)dy$$

$$\le C\int_\Omega G_\mu(x,z)^s\delta(z)^{\alpha(1-s)}\int_\Omega \delta(y)^\alpha\frac{G_\mu(x,y)}{\delta(x)}\left(\frac{G_\mu(y,z)}{\delta(z)^\alpha}\right)^{p-s}dyd\tau(z) \tag{3.11}$$

$$+ C\int_\Omega G_\mu(x,z)^s\delta(z)^{\alpha(1-s)}\int_\Omega \delta(y)^\alpha\left(\frac{G_\mu(x,y)}{\delta(x)}\right)^{1-s}\left(\frac{G_\mu(y,z)}{\delta(z)^\alpha}\right)^{p}dyd\tau(z) \tag{3.12}$$

$$\le C\int_\Omega G_\mu(x,z)^s\delta(z)^{\alpha(1-s)}\int_\Omega \delta(y)^\alpha\left(\frac{G_\mu(x,y)}{\delta(x)}\right)^{p-s+1}dyd\tau(z) \tag{3.13}$$

$$+ \int_\Omega G_\mu(x,z)^s\delta(z)^{\alpha(1-s)}\int_\Omega \delta(y)^\alpha\left(\frac{G_\mu(y,z)}{\delta(z)^\alpha}\right)^{p-s+1}dyd\tau(z) \tag{3.14}$$

$$\le C\int_\Omega \left(\frac{G_\mu(x,z)}{\delta(z)^\alpha}\right)^{s}\delta(z)^\alpha d\tau(z) \tag{3.15}$$

$$\le C\left(\int_\Omega G_\mu(x,z)d\tau(z)\right)^s. \tag{3.16}$$

Here (3.11) and (3.12) follow from (3.1), (3.13) and (3.14) follow from Hölder inequality, and (3.15) follows from Lemma 2.2, Hölder inequality and (3.9).

Note that the above approach can be applied to the case $p = 1$.

If $0 \le p < 1$ then

$$\mathbb{G}_\mu[\mathbb{G}_\mu[\tau]^p] \le C(\mathbb{G}_\mu[1] + \mathbb{G}_\mu[\mathbb{G}_\mu[\tau]]) \le C(\mathbb{G}_\mu[1] + \mathbb{G}_\mu[\tau]^s)$$

Then (3.10) follows by a similar argument as in the proof of Lemma 3.3. □

3.3. **Capacities and existence results.** For $a \ge 0$, $0 \le \theta \le \beta < N$ and $s > 1$, let $N_{\theta,\beta}$, $\mathbb{N}_{\theta,\beta}$ and $\mathrm{Cap}^a_{\mathbb{N}_{\theta,\beta},s}$ be defined as in (1.22), (1.23) and (1.24) respectively.

In this section we recall some results in [5, Section 2].

We recall below the definition of the capacity associated to $\mathbb{N}_{\theta,\beta}$ (see [11]).

**Definition 3.8.** *Let $a \ge 0$, $0 \le \theta \le \beta < N$ and $s > 1$. Define $Cap^a_{\mathbb{N}_{\theta,\beta},s}$ by*

$$\mathrm{Cap}^a_{\mathbb{N}_{\theta,\beta},s}(E) := \inf\left\{\int_{\overline{\Omega}}\delta^a\phi^s dy : \ \ \phi \ge 0, \ \ \mathbb{N}_{\theta,\beta}[\delta^a\phi] \ge \chi_E\right\},$$

*for any Borel set $E \subset \overline{\Omega}$.*

Clearly we have

$$\mathrm{Cap}^a_{\mathbb{N}_{\theta,\beta},s}(E) = \inf\left\{\int_{\overline{\Omega}}\delta^{-a(s-1)}\phi^s dy : \ \ \phi \ge 0, \ \ \mathbb{N}_{\theta,\beta}[\phi] \ge \chi_E\right\},$$

for any Borel set $E \subset \overline{\Omega}$. Furthermore we have by [1, Theorem 2.5.1]

$$\left(\mathrm{Cap}^a_{\mathbb{N}_{\theta,\beta},s}(E)\right)^{\frac{1}{s}} = \inf\left\{\omega(E) : \ \omega \in \mathfrak{M}^+_b(\overline{\Omega}), \ ||\mathbb{N}_{\theta,\beta}[\omega]||_{L^{s'}(\overline{\Omega};\delta^a)} \le 1\right\}, \tag{3.17}$$

for any compact set $E \subset \overline{\Omega}$ where $s'$ is the conjugate exponent of $s$.

Using [5, Theorem 2.6], we obtain easily the following result.

**Proposition 3.9.** *Let $p > 1$, $\sigma > 0$ and $\tau \in \mathfrak{M}^+(\Omega;\delta^\alpha)$. Then the following statements are equivalent.*

1. *There exists $C > 0$ such that the following inequality hold*

$$\int_E \delta^\alpha d\tau \leq C \mathrm{Cap}^{(p+1)\alpha}_{\mathbb{N}_{2\alpha,2},p'}(E),$$

*for any Borel $E \subset \overline{\Omega}$.*

2. *There exists a constant $C > 0$ such that* (1.19) *holds.*
3. *Problem* (3.5) *has a positive weak solution for $\sigma > 0$ small enough.*

*Proof.* First we note that

$$G(x,y) = \delta(x)^\alpha \delta(y)^\alpha N_{2\alpha,2}(x,y), \quad \forall x,y \in \Omega,\, x \neq y.$$

Thus the inequality

$$\mathbb{G}_\mu[\mathbb{G}_\mu[\tau]^p] \leq C\mathbb{G}_\mu[\tau] \quad \text{ a.e. in } \Omega$$

is equivalent to

$$\mathbb{N}_{2\alpha,2}[\delta^{(p+1)\alpha}(y)\mathbb{N}_{2\alpha,2}[\tilde{\tau}]^p(y)](x) \leq C\mathbb{N}_{2\alpha,2}[\tilde{\tau}](x) \quad \text{ a.e. in } \Omega,$$

where $d\tilde{\tau}(y) = \delta^\alpha(y)d\tau(y)$.

Now notice that if $u$ is a positive solution of (3.5) then by Proposition 2.12 we have that $u = \mathbb{G}_\mu[u^p] + \mathbb{G}_\mu[\tau]$ which implies that

$$\frac{u}{\delta(x)^\alpha} \approx \mathbb{N}_{2\alpha,2}\left[\delta^{(p+1)\alpha}\left(\frac{u}{\delta^\alpha}\right)^p\right](x) + \sigma\mathbb{N}_{2\alpha,2}[\tilde{\tau}](x),$$

the desired results follow by [5, Theorem 2.6] and [5, Proposition 2.7]. □

Let us now give a result which implies the existence for the problem (3.5).

**Lemma 3.10.** *Let $1 < p < p_\mu$. Then*

$$\inf_{\xi\in\Omega} \mathrm{Cap}^{(p+1)\alpha}_{\mathbb{N}_{2\alpha,2},p'}(\{\xi\}) > 0.$$

*Proof.* By (3.17) it is enough to show that

$$\sup_{\xi\in\Omega} ||\mathbb{N}_{2\alpha,2}[\delta_\xi]||_{L^p(\Omega;\delta^{(p+1)\alpha})} < C < \infty,$$

which is equivalent to

$$\sup_{\xi\in\Omega}\left|\left|\frac{G_\mu(\cdot,\xi)}{\delta(\xi)^\alpha}\right|\right|_{L^p(\Omega;\delta^\alpha)} < C. \tag{3.18}$$

The result follows by Lemma 2.2 and (2.1). □

3.4. **Boundary value problem.** Estimate (1.18) is a necessary and sufficient condition for the existence of weak solutions of

$$\begin{cases} -L_\mu u = u^p & \text{in } \Omega, \\ \operatorname{tr}(u) = \varrho\nu & \text{on } \partial\Omega. \end{cases} \tag{3.19}$$

**Proposition 3.11.** [5, Theorem 4.1] *Let $p > 1$, $\varrho > 0$ and $\nu \in \mathfrak{M}^+(\partial\Omega)$. Then, the following statements are equivalent.*

*1. There exists $C > 0$ such that the following inequality holds*

$$\nu(F) \leq C\mathrm{Cap}^{\partial\Omega}_{1-\alpha+\frac{1+\alpha}{p},p'}(F)$$

*for any Borel $F \subset \partial\Omega$.*

*2. There exists $C > 0$ such that (1.18) holds.*

*3. Problem (3.19) has a positive weak solution for $\varrho > 0$ small enough.*

**Lemma 3.12.** *Let $\nu \in \mathfrak{M}^+(\partial\Omega)$ and $0 < p < p_\mu$. Then there exists a constant $C > 0$ such that (1.18) holds.*

*Proof.* We first assume that $1 < p < p_\mu$. Let $\xi \in \partial\Omega$; we have $\delta_\xi(F) < c\mathrm{Cap}^{\partial\Omega}_{1-\alpha+\frac{1+\alpha}{p},p'}(F)$ for every $F \subset \partial\Omega$ where $c$ is independent of $\xi$. By Proposition 3.11, (1.18) holds with $\nu$ replaced by $\delta_\xi$ and with the constant $C$ independent of $\xi$. By taking integral over $\xi \in \partial\Omega$, we get (1.18).

Next, if $p \in (0, 1]$, we choose $s > 1$ such that $1 < ps < p_\mu$. By Young's inequality,

$$\mathbb{G}_\mu[\mathbb{K}_\mu[\nu]^p] \leq C(\mathbb{G}_\mu[1] + \mathbb{G}_\mu[\mathbb{K}_\mu[\nu]^{ps}]) \leq C(\mathbb{G}_\mu[1] + \mathbb{K}_\mu[\nu]). \tag{3.20}$$

This, combined with the inequality $\mathbb{G}_\mu[1] \leq c\delta^\alpha \leq c'\mathbb{K}_\mu[\nu]$ a.e. in $\Omega$ leads to (1.18). □

**Proposition 3.13.** *Let $p > 0$, $\varrho > 0$ and $\nu \in \mathfrak{M}^+(\partial\Omega)$.*

(i) *Assume there exists a constant $C > 0$ such that (1.18) holds. Then problem (3.19) admits a weak solution $u$ satisfying*

$$\mathbb{K}_\mu[\varrho\nu] \leq u \leq C\mathbb{K}_\mu[\varrho\nu] \quad \text{a.e. in } \Omega, \tag{3.21}$$

*with another constant $C > 0$, for any $\varrho > 0$ small enough if $p > 1$, for any $\varrho > 0$ if $p \in (0, 1)$.*

(ii) *Assume $p > 1$ and problem (3.19) admits a weak solution. Then (1.18) holds with $C = \frac{1}{p-1}$.*

(iii) *Assume $0 < p < p_\mu$. Then there exists a constant $C > 0$ such that for any weak solution $u$ of (3.19) there holds*

$$\mathbb{K}_\mu[\varrho\nu] \leq u \leq C(\mathbb{K}_\mu[\varrho\nu] + \delta^\alpha) \quad \text{a.e. in } \Omega. \tag{3.22}$$

*Proof.* By using an argument as in the proof of Proposition 3.4, Proposition 3.5 and Proposition 3.6, we obtained the desired results. □

The above results allow to study elliptic equations with interior and boundary measures.

**Proposition 3.14.** *Let $p > 0$, $\sigma > 0$, $\varrho > 0$ and $\tau \in \mathfrak{M}^+(\Omega; \delta^\alpha)$ and $\nu \in \mathfrak{M}^+(\partial\Omega)$. If (1.18) and (1.19) hold then problem (1.17) admits a weak solution $u$ satisfying (1.20) for $\sigma > 0$ and $\varrho > 0$ small enough if $p > 1$, for any $\sigma > 0$ and $\varrho > 0$ if $0 < p < 1$.*

*Furthermore if $0 < p < p_\mu$ there exists a constant $C > 0$ such that for any weak solution $u$ of (1.17) estimate (1.21) holds.*

*Proof.* We adapt the argument in the proof of [4, Theorem 3.13]. Put $v := u - \mathbb{K}_\mu[\varrho\nu]$ then $v$ satisfies

$$\begin{cases} -L_\mu v = (v + \mathbb{K}_\mu[\varrho\nu])^p + \sigma\tau \quad \text{in } \Omega, \\ \operatorname{tr}(v) = 0. \end{cases} \tag{3.23}$$

Consider the following problem

$$\begin{cases} -L_\mu w = c_p w^p + c_p(\mathbb{K}_\mu[\varrho\nu])^p + \sigma\tau \quad \text{in } \Omega, \\ \operatorname{tr}(w) = 0 \end{cases} \tag{3.24}$$

where $c_p := \max\{1, 2^{p-1}\}$. Since (1.18) holds, it follows that $\mathbb{K}_\mu[\nu]^p \in L^1(\Omega;\delta^\alpha)$. Since (1.19) holds, we infer from Proposition 3.4 that problem (3.24) admits a weak solution $w$ for $\sigma > 0$ and $\varrho > 0$ small enough if $p > 1$, for any $\sigma > 0$ and $\varrho > 0$ if $0 < p < 1$. Notice that $w$ is a supersolution of (3.24), we infer that there is a weak solution $v$ of (3.23) satisfying $v \le w$ a. e. in $\Omega$. By Proposition 3.4 and (1.18), we get

$$w \le c\mathbb{G}_\mu[\mathbb{K}_\mu[\varrho\nu]^p + \sigma\tau] \le c'(\mathbb{G}_\mu[\sigma\tau] + \mathbb{K}_\mu[\varrho\nu]) \quad \text{a.e. in } \Omega.$$

This implies (1.20).

If $0 < p < p_\mu$ then (1.21) follows from Proposition 3.6 and Proposition 3.13 (iii). □

**Proof of Theorem B.** Statements (i) and (ii) follow from Lemma 3.12 and Lemma 3.3 respectively. Statement (iii) follows from Proposition 3.14. Statement (iv) follows from Proposition 3.5 and Proposition 3.13 (ii). Statement (v) is derived from Proposition 3.14 (ii). □

**Proof of Theorem C.** The implications (i) $\Longleftrightarrow$ (ii) $\Longrightarrow$ (iii) follow from Proposition 3.11, Proposition 3.9 and Proposition 3.14. We will show that (iii) $\Longrightarrow$ (ii). Since (1.17) has a weak solution for $\sigma > 0$ small and $\varrho > 0$ small, it follows that (3.5) admits a solution for $\sigma > 0$ small and (3.19) admits a solution for $\varrho > 0$ small. Due to Proposition 3.11 and Proposition 3.9, we derive (1.19) and (1.18). This completes the proof. □

## 4. Elliptic systems: the power case

Let $\mu \in (0, \frac{1}{4}]$. In this section, we deal with system (1.32). We recall that $p_\mu$ is defined in (1.31) and

$$q := \tilde{p}\frac{p+1}{\tilde{p}+1}, \quad \tilde{q} := p\frac{\tilde{p}+1}{p+1}.$$

Without loss of generality, we can assume that $0 < p \le \tilde{p}$. Then $p \le q \le \tilde{q} \le \tilde{p}$ if $p\tilde{p} \ge 1$. Put

$$t_\mu := \tilde{p}\,(p - p_\mu + 1)\,.$$

Notice that if $q < p_\mu$ then $t_\mu < q < p_\mu$.

**Lemma 4.1.** *Let $p > 0$, $\tilde{p} > 0$ and $\tau \in \mathfrak{M}^+(\Omega;\delta^\alpha)$. Assume $q < p_\mu$. Then for any $t \in (\max(0,t_\mu), \tilde{p}]$, there exists a positive constant $c = c(N,p,\tilde{p},\mu,t,\tau)$ (independent of $\tau$ if $p > 1$) such that*

$$\mathbb{G}_\mu[\mathbb{G}_\mu[\tau]^p]^{\tilde{p}} \le c\mathbb{G}_\mu[\tau]^t. \tag{4.1}$$

*In particular,*

$$\mathbb{G}_\mu[\mathbb{G}_\mu[\tau]^p]^{\tilde{p}} \le C\mathbb{G}_\mu[\tau]^q, \tag{4.2}$$

$$\mathbb{G}_\mu[\mathbb{G}_\mu[\mathbb{G}_\mu[\tau]^p]^{\tilde{p}}] \le C\mathbb{G}_\mu[\tau] \tag{4.3}$$

*where* $C = C(N, p, \tilde{p}, \mu, \tau)$.

*Proof.* Since $q < p_\mu$, it follows that $p < p_\mu$, hence $\max(0, p - p_\mu + 1) < 1$. Let $t \in (\max(0, t_\mu, \tilde{p}]$ then $\max(0, p - p_\mu + 1) < \frac{t}{\tilde{p}} \le 1$. By applying Lemma 3.7 with $s$ replaced by $\frac{t}{\tilde{p}}$ respectively in order to obtain

$$\mathbb{G}_\mu[\mathbb{G}_\mu[\tau]^p] \le c\mathbb{G}_\mu[\tau]^{\frac{t}{\tilde{p}}},$$

which implies (4.1). Since $t_\mu < q \le \tilde{p}$, by taking $t = q$ in (4.1) we obtain (4.2). Next, since $q < p_\mu$, by apply Lemma 3.7 with $\gamma$ replaced by $\tilde{\gamma}$ and (4.2), we get

$$\mathbb{G}_\mu[\mathbb{G}_\mu[\mathbb{G}_\mu[\tau]^p]^{\tilde{p}}] \le C\mathbb{G}_\mu[\mathbb{G}_\mu[\tau]^q] \le C\mathbb{G}_\mu[\tau].$$

□

**Lemma 4.2.** *Let* $p > 0$, $\tilde{p} > 0$, $\tau, \tilde{\tau} \in \mathfrak{M}^+(\Omega; \delta^\alpha)$ *and* $\nu, \tilde{\nu} \in \mathfrak{M}^+(\partial\Omega)$. *Assume that there exist positive functions* $U \in L^{\tilde{p}}(\Omega; \delta^\alpha)$ *and* $V \in L^p(\Omega; \delta^\alpha)$ *such that*

$$\begin{aligned} U &\ge \mathbb{G}_\mu[(V + \mathbb{K}_\mu[\tilde{\varrho}\tilde{\nu}])^p] + \mathbb{G}_\mu[\sigma\tau], \\ V &\ge \mathbb{G}_\mu[(U + \mathbb{K}_\mu[\varrho\nu])^{\tilde{p}}] + \mathbb{G}_\mu[\tilde{\sigma}\tilde{\tau}] \end{aligned} \tag{4.4}$$

*in* $\Omega$. *Then there exists a weak solution* $(u, v)$ *of* (1.32) *such that*

$$\begin{aligned} \mathbb{G}_\mu[\sigma\tau] + \mathbb{K}_\mu[\varrho\nu] &\le u \le U, \\ \mathbb{G}_\mu[\tilde{\sigma}\tilde{\tau}] + \mathbb{K}_\mu[\tilde{\varrho}\tilde{\nu}] &\le v \le V. \end{aligned} \tag{4.5}$$

*Proof.* Put $u_0 := 0$ and

$$\begin{cases} v_{n+1} := \mathbb{G}_\mu[u_n^{\tilde{p}}] + \mathbb{G}_\mu[\tilde{\sigma}\tilde{\tau}] + \mathbb{K}_\mu[\tilde{\varrho}\tilde{\nu}], & n \ge 0, \\ u_n := \mathbb{G}_\mu[v_n^p] + \mathbb{G}_\mu[\sigma\tau] + \mathbb{K}_\mu[\varrho\nu], & n \ge 1. \end{cases} \tag{4.6}$$

We see that $0 \le v_1 = \mathbb{G}_\mu[\tilde{\sigma}\tilde{\tau}] + \mathbb{K}_\mu[\tilde{\varrho}\tilde{\nu}] \le V$. It is easy to see that $\{u_n\}$ and $\{v_n\}$ are nondecreasing sequences, $0 \le u_n \le U$ and $0 \le v_n \le V$ in $\Omega$. By monotone convergence theorem, there exist $u \in L^{\tilde{p}}(\Omega; \delta^\alpha)$ and $v \in L^p(\Omega; \delta^\alpha)$ such that $u_n \to u$ in $L^1(\Omega)$, $v_n \to v$ in $L^1(\Omega)$, $u_n^{\tilde{p}} \to u^{\tilde{p}}$ in $L^1(\Omega; \delta^\alpha)$, $v_n^p \to v^p$ in $L^1(\Omega; \delta^\alpha)$. Moreover $u \le U$ and $v \le V$ in $\Omega$. By letting $n \to \infty$ in (4.6), we obtain

$$\begin{cases} v = \mathbb{G}_\mu[u^{\tilde{p}}] + \mathbb{G}_\mu[\tilde{\sigma}\tilde{\tau}] + \mathbb{K}_\mu[\tilde{\varrho}\tilde{\nu}], \\ u = \mathbb{G}_\mu[v^p] + \mathbb{G}_\mu[\sigma\tau] + \mathbb{K}_\mu[\varrho\nu]. \end{cases} \tag{4.7}$$

Thus $(u, v)$ is a weak solution of (1.32) and satisfies (4.5). □

**Proof of Theorem E.** We first show that the following system has weak a solution

$$\begin{cases} -L_\mu w = (\tilde{w} + \mathbb{K}_\mu[\tilde{\varrho}\tilde{\nu}])^p + \sigma\tau & \text{in } \Omega, \\ -L_\mu \tilde{w} = (w + \mathbb{K}_\mu[\varrho\nu])^{\tilde{p}} + \tilde{\sigma}\tilde{\tau} & \text{in } \Omega, \\ \operatorname{tr}(u) = \operatorname{tr}(v) = 0. \end{cases} \tag{4.8}$$

Fix $\vartheta_i > 0$, $(i = 1, 2, 3, 4)$ and set

$$\Psi := \mathbb{G}_\mu[\vartheta_1\tau + \mathbb{K}_\mu[\vartheta_2\tilde{\nu}]^p]^{\tilde{p}} + \mathbb{K}_\mu[\vartheta_3\nu]^{\tilde{p}} + \vartheta_4\tilde{\tau}.$$

For $\kappa \in (0, 1]$, put

$$\sigma := \kappa^{\frac{1}{\tilde{p}}}\vartheta_1, \quad \tilde{\sigma} := \kappa\vartheta_4, \quad \varrho := \kappa^{\frac{1}{\tilde{p}}}\vartheta_3, \quad \tilde{\varrho} := \kappa^{\frac{1}{p\tilde{p}}}\vartheta_2.$$

Then from the assumption, we deduce that $\Psi \in \mathfrak{M}^+(\Omega;\delta^\alpha)$. By Lemma 4.1,

$$\mathbb{G}_\mu[\mathbb{G}_\mu[\mathbb{G}_\mu[\Psi]^p]^{\tilde p}] \leq C\mathbb{G}_\mu[\Psi] \tag{4.9}$$

where $C = C(N,p,\tilde p,\mu,\sigma,\tilde\sigma,\kappa,\tau,\tilde\tau)$. Set

$$V := A\mathbb{G}_\mu[\kappa\Psi] \quad \text{and} \quad U := \mathbb{G}_\mu[(V+\mathbb{K}_\mu[\tilde\varrho\tilde\nu])^p + \sigma\tau]$$

where $A>0$ will be determined later on. We have

$$\begin{aligned}
&(U+\mathbb{K}_\mu[\varrho\nu])^{\tilde p} + \tilde\sigma\tilde\tau \\
&\leq c\left\{\mathbb{G}_\mu\Big[(a_3\kappa\mathbb{G}_\mu[\Psi]+\mathbb{K}_\mu[\tilde\varrho\tilde\nu])^p\Big]^{\tilde p} + \mathbb{G}_\mu[\sigma\tau]^{\tilde p} + \mathbb{K}_\mu[\varrho\nu]^{\tilde p}\right\} + \tilde\sigma\tilde\tau \\
&\leq c\left\{A^{p\tilde p}\kappa^{p\tilde p}\mathbb{G}_\mu[\mathbb{G}_\mu[\Psi]^p]^{\tilde p} + \mathbb{G}_\mu[\mathbb{K}_\mu[\tilde\varrho\tilde\nu]^p]^{\tilde p}\right\} + c\,\mathbb{G}_\mu[\sigma\tau]^{\tilde p} + c\,\mathbb{K}_\mu[\varrho\nu]^{\tilde p} + \tilde\sigma\tilde\tau
\end{aligned}$$

where $c = c(p,\tilde p)$. It follows that

$$\mathbb{G}_\mu[(U+\mathbb{K}_\mu[\varrho\nu])^{\tilde p} + \tilde\sigma\tilde\tau] \leq I_1 + I_2 \tag{4.10}$$

where

$$\begin{aligned}
I_1 &:= c\,A^{p\tilde p}\kappa^{p\tilde p}\mathbb{G}_\mu\Big[\mathbb{G}_\mu[\mathbb{G}_\mu[\Psi]^p]^{\tilde p}\Big] + c\,\mathbb{G}_\mu\Big[\mathbb{G}_\mu[\mathbb{K}_\mu[\tilde\varrho\tilde\nu]^p]^{\tilde p}\Big], \\
I_2 &:= c\,\mathbb{G}_\mu[\mathbb{G}_\mu[\sigma\tau]^{\tilde p}] + c\,\mathbb{G}_\mu[\mathbb{K}_\mu[\varrho\nu]^{\tilde p}] + \mathbb{G}_\mu[\tilde\sigma\tilde\tau].
\end{aligned}$$

We first estimate $I_1$. Observe that

$$\mathbb{G}_\mu\Big[\mathbb{G}_\mu[\mathbb{K}_\mu[\tilde\varrho\tilde\nu]^p]^{\tilde p}\Big] = \mathbb{G}_\mu\Big[\mathbb{G}_\mu[\mathbb{K}_\mu[\vartheta_2\tilde\nu]^p]^{\tilde p}\Big] \leq \mathbb{G}_\mu[\kappa\Psi].$$

This, together with (4.9) implies

$$I_1 \leq c\,(A^{p\tilde p}\kappa^{p\tilde p-1}+1)\mathbb{G}_\mu[\kappa\Psi]. \tag{4.11}$$

Next it is easy to see that

$$I_2 \leq c\,\mathbb{G}_\mu[\kappa\Psi]. \tag{4.12}$$

By collecting (4.10), (4.11) and (4.12), we obtain

$$\mathbb{G}_\mu[(U+\mathbb{K}_\mu[\varrho\nu])^{\tilde p} + \tilde\sigma\tilde\tau] \leq c(A^{p\tilde p}\kappa^{p\tilde p-1}+1)\mathbb{G}_\mu[\kappa\Psi] \tag{4.13}$$

with another constant $c$. We will choose $A$ and $\kappa$ such that

$$c\,(A^{p\tilde p}\kappa^{p\tilde p-1}+1) \leq A. \tag{4.14}$$

If $p\tilde p > 1$ then we can choose $A>0$ large enough and then choose $\kappa>0$ small enough (depending on $A$) such that (4.14) holds. If $p\tilde p<1$ then for any $\kappa>0$ there exists $A$ large enough such that (4.14) holds. For such $A$ and $\kappa$, we obtain

$$\mathbb{G}_\mu[(U+\mathbb{K}_\mu[\varrho\nu])^{\tilde p} + \tilde\sigma\tilde\tau] \leq V.$$

By Lemma 4.2, there exists a weak solution $(w,\tilde w)$ of (4.8) for $\sigma>0$, $\tilde\sigma>0$, $\nu>0$, $\tilde\nu>0$ small if $p\tilde p>1$, for any $\sigma>0$, $\tilde\sigma>0$, $\nu>0$, $\tilde\nu>0$ if $p\tilde p<1$. Moreover, $(w,\tilde w)$ satisfies

$$\tilde w \approx \mathbb{G}_\mu[\omega], \tag{4.15}$$

$$w \approx \mathbb{G}_\mu[(\mathbb{G}_\mu[\omega]+\mathbb{K}_\mu[\tilde\nu])^p] + \mathbb{G}_\mu[\tau] \tag{4.16}$$

where $C = C(N,p,\tilde p,\mu,\Omega,\sigma,\tilde\sigma,\tau,\tilde\tau)$.

Next put $u := w + \mathbb{K}_\mu[\varrho\nu]$ and $v := \tilde w + \mathbb{K}_\mu[\tilde\varrho\tilde\nu]$ then $(u,v)$ is a weak solution of (1.32). Moreover (1.33) and (1.34) follow directly from (4.15) and (4.16). □

**Proof of Theorem F.** Put $\tau^* := \max\{\tau, \tilde{\tau}\}$ and $\nu^* := \max\{\nu, \tilde{\nu}\}$. Fix $\vartheta > 0$, $\tilde{\vartheta} > 0$ and for $\kappa \in (0,1]$, put $\sigma = \varrho = (\kappa\vartheta)^{\frac{1}{\tilde{p}}}$ and $\tilde{\sigma} = \kappa\tilde{\vartheta}$, $\tilde{\varrho} = (\kappa\tilde{\vartheta})^{\frac{1}{p\tilde{p}}}$. Set

$$\tau^{\#} := \vartheta\tau + \tilde{\vartheta}\tilde{\tau} \quad \text{and} \quad \nu^{\#} := \vartheta\nu + \tilde{\vartheta}\tilde{\nu}$$

then $\tau^{\#} \le (\vartheta + \tilde{\vartheta})\tau^*$ and $\nu^{\#} \le (\vartheta + \tilde{\vartheta})\nu^*$.

Put $V := A(\mathbb{G}_\mu[\kappa\tau^{\#}] + \mathbb{K}_\mu[\kappa\nu^{\#}])$ where $A > 0$ will be determined later on and put $U := \mathbb{G}_\mu[(V + \mathbb{K}_\mu[\tilde{\varrho}\tilde{\nu}])^p + \sigma\tau]$.

We have

$$U^{\tilde{p}} + \tilde{\sigma}\tilde{\tau} \le c\,A^{p\tilde{p}}\kappa^{p\tilde{p}}\Big\{\mathbb{G}_\mu[\mathbb{G}_\mu[\tau^{\#}]^p]^{\tilde{p}} + \mathbb{G}_\mu[\mathbb{K}_\mu[\nu^{\#}]^p]^{\tilde{p}}\Big\} + c\,\sigma^{\tilde{p}}\mathbb{G}_\mu[\tau]^{\tilde{p}} + c\,\varrho^{\tilde{p}}\mathbb{K}_\mu[\nu]^{\tilde{p}} + c\,\tilde{\sigma}\tilde{\tau}.$$

with $c = c(p, \tilde{p})$. It follows that

$$\mathbb{G}_\mu[(U + \mathbb{K}_\mu[\varrho\nu])^{\tilde{p}} + \tilde{\sigma}\tilde{\tau}] \le c\,(J_1 + J_2) \tag{4.17}$$

where

$$J_1 := A^{p\tilde{p}}\kappa^{p\tilde{p}}\Big\{\mathbb{G}_\mu[\mathbb{G}_\mu[\mathbb{G}_\mu[\tau^{\#}]^p]^{\tilde{p}}] + \mathbb{G}_\mu[\mathbb{G}_\mu[\mathbb{K}_\mu[\nu^{\#}]^p]^{\tilde{p}}]\Big\},$$

$$J_2 := \sigma^{\tilde{p}}\mathbb{G}_\mu[\mathbb{G}_\mu[\tau]^{\tilde{p}}] + \varrho^{\tilde{p}}\mathbb{G}_\mu[\mathbb{K}_\mu[\nu]^{\tilde{p}}] + \tilde{\sigma}\mathbb{G}_\mu[\tilde{\tau}] + \tilde{\varrho}\mathbb{K}_\mu[\tilde{\nu}].$$

We first estimate $J_1$. We have

$$J_1 \le A^{p\tilde{p}}\kappa^{p\tilde{p}}(\vartheta + \tilde{\vartheta})^{p\tilde{p}}\Big\{\mathbb{G}_\mu[\mathbb{G}_\mu[\mathbb{G}_\mu[\tau^*]^p]^{\tilde{p}}] + \mathbb{G}_\mu[\mathbb{G}_\mu[\mathbb{K}_\mu[\nu^*]^p]^{\tilde{p}}]\Big\}.$$

By (1.35), (1.36) and Proposition 3.11, Proposition 3.9 we infer that

$$J_1 \le c\,A^{p\tilde{p}}\kappa^{p\tilde{p}}(\vartheta + \tilde{\vartheta})^{p\tilde{p}}(\mathbb{G}_\mu[\tau^*] + \mathbb{K}_\mu[\nu^*])$$

where $c$ is a positive constant. Therefore

$$J_1 \le c\,A^{p\tilde{p}}\kappa^{p\tilde{p}}(\vartheta + \tilde{\vartheta})^{p\tilde{p}}\max(\vartheta^{-1}, \tilde{\vartheta}^{-1})(\mathbb{G}_\mu[\tau^{\#}] + \mathbb{K}_\mu[\nu^{\#}]). \tag{4.18}$$

We next estimate $J_2$. Again by (1.35), (1.36) and Proposition 3.11, Proposition 3.9, we deduce

$$\begin{aligned} J_2 &\le c\,(\sigma^{\tilde{p}}\mathbb{G}_\mu[\tau] + \varrho^{\tilde{p}}\mathbb{K}_\mu[\nu] + \tilde{\sigma}\mathbb{G}_\mu[\tilde{\tau}] + \tilde{\varrho}\mathbb{K}_\mu[\tilde{\nu}]) \\ &= c\,\kappa(\mathbb{G}_\mu[\tau^{\#}] + \mathbb{K}_\mu[\nu^{\#}]). \end{aligned} \tag{4.19}$$

Combining (4.17), (4.18) and (4.19) implies

$$\mathbb{G}_\mu[U^{\tilde{p}} + \tilde{\sigma}\tilde{\tau}] + \mathbb{K}_\mu[\tilde{\varrho}\tilde{\nu}] \le C(A^{p\tilde{p}}\kappa^{p\tilde{p}-1} + 1)(\mathbb{G}_\mu[\kappa\tau^{\#}] + \mathbb{K}_\mu[\kappa\nu^{\#}]) \tag{4.20}$$

where $C$ is another positive constant. We choose $A > 0$ and $\kappa > 0$ such that

$$C(A^{p\tilde{p}}\kappa^{p\tilde{p}-1} + 1) \le A. \tag{4.21}$$

Since $p\tilde{p} > 1$, one can choose $A$ large enough and then choose $\kappa > 0$ small enough such that (4.21) holds. For such $A$ and $\kappa$, we have

$$\mathbb{G}_\mu[U^{\tilde{p}} + \tilde{\sigma}\tilde{\tau}] + \mathbb{K}_\mu[\tilde{\varrho}\tilde{\nu}] \le V.$$

By Lemma 4.2, there exists a weak solution $(u, v)$ of (1.32) which satisfies (1.37). $\square$

## 5. General nonlinearities

5.1. **Absorption case.** In this section we treat system (1.38) with $\epsilon = -1$. We recall that $\Lambda_g$ and $\Lambda_{\tilde g}$ are defined in (1.39).

**Proof of Theorem G.**
**Step 1**: We claim that

$$\int_\Omega g(\mathbb{K}_\mu[|\tilde\nu|] + \mathbb{G}_\mu[|\tilde\sigma|])\delta^\alpha dx + \int_\Omega \tilde g(\mathbb{K}_\mu[|\nu|] + \mathbb{G}_\mu[|\tau|])\delta^\alpha dx < \infty. \tag{5.1}$$

For $\lambda > 0$, set $\tilde A_\lambda := \{x \in \Omega : \mathbb{K}_\mu[|\tilde\nu|] + \mathbb{G}_\mu[|\tilde\tau|] > \lambda\}$ and $a(\lambda) := \int_{\tilde A_\lambda} \delta^\alpha dx$. We write

$$\begin{aligned}\|g(\mathbb{K}_\mu[|\tilde\nu|] + \mathbb{G}_\mu[|\tilde\tau|])\|_{L^1(\Omega;\delta^\alpha)} &= \int_{\tilde A_1} g(\mathbb{K}_\mu[|\tilde\nu|] + \mathbb{G}_\mu[|\tilde\tau|])\delta^\alpha dx + \int_{\tilde A_1^c} g(\mathbb{K}_\mu[|\tilde\nu|] + \mathbb{G}_\mu[|\tilde\tau|])\delta^\alpha dx \\ &\le \int_{\tilde A_1} g(\mathbb{K}_\mu[|\tilde\nu|] + \mathbb{G}_\mu[|\tilde\tau|])\delta^\alpha dx + g(1)\int_\Omega \delta^\alpha dx.\end{aligned} \tag{5.2}$$

We have

$$\int_{\tilde A_1} g(\mathbb{K}_\mu[|\tilde\nu|] + \mathbb{G}_\mu[|\tilde\tau|])\delta^\alpha dx = a(1)g(1) + \int_1^\infty a(s)dg(s).$$

On the other hand, by (2.2) and Proposition 2.4 one gets, for every $s > 0$,

$$a(s) \le C\left(\|\mathbb{K}_\mu[|\tilde\nu|]\|^{p_\mu}_{L^{p_\mu}_w(\Omega;\delta^\alpha)} + \|\mathbb{G}_\mu[|\tilde\tau|]\|^{p_\mu}_{L^{p_\mu}_w(\Omega;\delta^\alpha)}\right) s^{-p_\mu} \le C\, s^{-p_\mu} \tag{5.3}$$

where $C = C(N, \mu, \Omega, \gamma, \|\tilde\nu\|_{\mathfrak{M}(\partial\Omega)}, \|\tilde\tau\|_{\mathfrak{M}(\Omega;\delta^\alpha)})$. Thus

$$a(1)g(1) + \int_1^\infty a(s)dg(s) \le C + C\int_1^\infty s^{-1-p_\mu} g(s)ds \le Cp_\mu\Lambda_g. \tag{5.4}$$

By combining the above estimates we obtain

$$\|g(\mathbb{K}_\mu[|\tilde\nu|] + \mathbb{G}_\mu[|\tilde\tau|])\|_{L^1(\Omega;\delta^\alpha)} \le Cp_\mu\Lambda_g + g(1)\int_\Omega \delta^\alpha dx \le C.$$

Similarly,

$$\|\tilde g(\mathbb{K}_\mu[|\nu|] + \mathbb{G}_\mu[|\tau|])\|_{L^1(\Omega;\delta^\alpha)} \le \tilde Cp_\mu\Lambda_{\tilde g} + \tilde g(1)\int_\Omega \delta^\alpha dx \le \tilde C.$$

Thus (5.1) follows directly.

**Step 2**: Existence.

Put $u_0 := \mathbb{K}_\mu[\nu] + \mathbb{G}_\mu[\tau]$. Let $v_0$ be the unique weak solution of the following problem

$$\begin{cases} -L_\mu v_0 + \tilde g(u_0) = \tilde\tau & \text{in } \Omega, \\ \operatorname{tr}(v_0) = \tilde\nu. \end{cases}$$

For any $k \ge 1$, since $g, \tilde g$ satisfy (5.1) there exist functions $u_k$ and $v_k$ satisfying

$$\begin{cases} -L_\mu u_k + g(v_{k-1}) = \tau & \text{in } \Omega, \\ -L_\mu v_k + \tilde g(u_k) = \tilde\tau & \text{in } \Omega, \\ \operatorname{tr}(u_k) = \nu, \quad \operatorname{tr}(v_k) = \tilde\nu. \end{cases} \tag{5.5}$$

Moreover

$$(5.6)\qquad \begin{aligned} u_k + \mathbb{G}_\mu[g(v_{k-1})] &= \mathbb{G}_\mu[\tau] + \mathbb{K}_\mu[\nu], \\ v_k + \mathbb{G}_\mu[\tilde{g}(u_k)] &= \mathbb{G}_\mu[\tilde{\tau}] + \mathbb{K}_\mu[\tilde{\nu}]. \end{aligned}$$

Since $g, \tilde{g} \geq 0$, it follows that, for every $k \geq 1$,

$$\mathbb{K}_\mu[\nu] + \mathbb{G}_\mu[\tau] - \mathbb{G}_\mu[g(\mathbb{K}_\mu[\tilde{\nu}] + \mathbb{G}_\mu[\tilde{\tau}])] \leq u_k \leq \mathbb{K}_\mu[\nu] + \mathbb{G}_\mu[\tau] = u_0$$

and

$$\mathbb{K}_\mu[\tilde{\nu}] + \mathbb{G}_\mu[\tilde{\tau}] - \mathbb{G}_\mu[\tilde{g}(\mathbb{K}_\mu[\nu] + \mathbb{G}_\mu[\tau])] \leq v_k \leq \mathbb{K}_\mu[\tilde{\nu}] + \mathbb{G}_\mu[\tilde{\tau}]$$

in $\Omega$. Now, suppose that for some $k \geq 1$, $u_k \leq u_{k-1}$. Since $g$ and $\tilde{g}$ are nondecreasing, we deduce that

$$(5.7)\qquad \begin{aligned} v_k &= \mathbb{K}_\mu[\tilde{\nu}] + \mathbb{G}_\mu[\tilde{\tau}] - \mathbb{G}_\mu[\tilde{g}(u_k)] \geq \mathbb{K}_\mu[\tilde{\nu}] + \mathbb{G}_\mu[\tilde{\tau}] - \mathbb{G}_\mu[\tilde{g}(u_{k-1})] = v_{k-1}, \\ u_{k+1} &= \mathbb{K}_\mu[\nu] + \mathbb{G}_\mu[\tau] - \mathbb{G}_\mu[g(v_k)] \leq \mathbb{K}_\mu[\nu] + \mathbb{G}_\mu[\tau] - \mathbb{G}_\mu[g(v_{k-1})] = u_k. \end{aligned}$$

This means that $\{v_k\}$ is nondecreasing and $\{u_k\}$ is nonincreasing. Hence, there exist $u$ and $v$ such that $u_k \downarrow u$ and $v_k \uparrow v$ in $\Omega$ and

$$\mathbb{K}_\mu[\nu] + \mathbb{G}_\mu[\tau] - \mathbb{G}_\mu[g(\mathbb{K}_\mu[\tilde{\nu}] + \mathbb{G}_\mu[\tilde{\tau}])] \leq u \leq \mathbb{K}_\mu[\nu] + \mathbb{G}_\mu[\tau],$$

$$\mathbb{K}_\mu[\tilde{\nu}] + \mathbb{G}_\mu[\tilde{\tau}] - \mathbb{G}_\mu[\tilde{g}(\mathbb{K}_\mu[\nu] + \mathbb{G}_\mu[\tau])] \leq v \leq \mathbb{K}_\mu[\tilde{\nu}] + \mathbb{G}_\mu[\tilde{\tau}].$$

Since $g$ and $\tilde{g}$ are continuous and nondecreasing, we infer from monotone convergence theorem and (5.1) that $g(v_k) \to g(v)$ in $L^1(\Omega; \delta^\alpha)$ and $\tilde{g}(u_k) \to \tilde{g}(u)$ in $L^1(\Omega; \delta^\alpha)$. As a consequence,

$$\begin{aligned} \mathbb{G}_\mu[\tilde{g}(u_k)] &\to \mathbb{G}_\mu[\tilde{g}(u)] \qquad \text{a.e. in } \Omega, \\ \mathbb{G}_\mu[g(v_k)] &\to \mathbb{G}_\mu[g(v)] \qquad \text{a.e. in } \Omega. \end{aligned}$$

By letting $k \to \infty$ in (5.6), we obtain the desired result. □

5.2. **Source case: subcriticality.** In this section for simplicity we consider system (1.38) with $\epsilon = 1$. Assume that $g(0) = \tilde{g}(0) = 0$. In preparation for proving Theorem H, we establish the following lemma:

**Lemma 5.1.** *Assume $\epsilon = 1$, $g$ and $\tilde{g}$ are bounded, nondecreasing and continuous functions in $\mathbb{R}$. Let $\tau$, $\tilde{\tau} \in \mathfrak{M}(\Omega; \delta^\alpha)$ and $\nu, \tilde{\nu} \in \mathfrak{M}(\partial\Omega)$. Assume there exist $a_1 > 0$, $b_1 > 0$ and $q_1 > 1$ such that* (1.40) *and* (1.41) *are satisfied. Then there exist $\lambda_*$, $\tilde{\lambda}_*$, $b_* > 0$ and $\varrho_* > 0$ depending on $N$, $\mu$, $\Omega$ $\gamma$, $\tilde{\gamma}$, $\Lambda_g$, $\Lambda_{\tilde{g}}$, $a_1$, $q_1$ such that the following holds. For every $b_1 \in (0, b_*)$ and $\sigma, \tilde{\sigma}, \tilde{\varrho}, \varrho \in (0, \varrho_*)$ the system*

$$(5.8)\qquad \begin{cases} -L_\mu u = g(v + \tilde{\varrho}\mathbb{K}_\mu[\tilde{\nu}] + \tilde{\sigma}\mathbb{G}_\mu[\tilde{\tau}]) & \textit{in } \Omega, \\ -L_\mu v = \tilde{g}(u + \varrho\mathbb{K}_\mu[\nu] + \sigma\mathbb{G}_\mu[\tau]) & \textit{in } \Omega, \\ \mathrm{tr}\,(u) = \mathrm{tr}\,(v) = 0 \end{cases}$$

*admits a weak solution $(u, v)$ satisfying*

$$(5.9)\qquad \begin{aligned} \|u\|_{L_w^{p_\mu}(\Omega;\delta^\alpha)} + \|u\|_{L^{q_1}(\Omega;\delta^{\alpha-1})} &\leq \lambda_*, \\ \|v\|_{L_w^{p_\mu}(\Omega;\delta^\alpha)} + \|v\|_{L^{q_1}(\Omega;\delta^{\alpha-1})} &\leq \tilde{\lambda}_*. \end{aligned}$$

*Proof.* Without loss of generality, we assume that $\|\tau\|_{\mathfrak{M}(\Omega;\delta^\alpha)} = \|\tilde\tau\|_{\mathfrak{M}(\Omega;\delta^\alpha)} = \|\nu\|_{\mathfrak{M}(\partial\Omega)} = \|\tilde\nu\|_{\mathfrak{M}(\partial\Omega)} = 1$. We shall use Schauder fixed point theorem to show the existence of positive weak solutions of (5.8). Define

$$\begin{aligned} &\mathbb{S}(w) := \mathbb{G}_\mu[g(w + \tilde\varrho\mathbb{K}_\mu[\tilde\nu] + \tilde\sigma\mathbb{G}_\mu[\tilde\tau])], \\ &\tilde{\mathbb{S}}(w) := \mathbb{G}_\mu[\tilde g(w + \varrho\mathbb{K}_\mu[\nu] + \sigma\mathbb{G}_\mu[\tau])], \quad \forall w \in L^1(\Omega). \end{aligned} \tag{5.10}$$

Set

$$\begin{aligned} &\mathbf{M}_1(w) := \|w\|_{L_w^{p_\mu}(\Omega;\delta^\alpha)}, \quad \forall w \in L_w^{p_\mu}(\Omega;\delta^\alpha), \\ &\tilde{\mathbf{M}}_1(w) := \|w\|_{L_w^{p_\mu}(\Omega;\delta^\alpha)}, \quad \forall w \in L_w^{p_\mu}(\Omega;\delta^\alpha), \\ &\mathbf{M}_2(w) := \|w\|_{L^{q_1}(\Omega;\delta^{\alpha-1})}, \quad \forall w \in L^{q_1}(\Omega;\delta^{\alpha-1}), \\ &\mathbf{M}(w) := \mathbf{M}_1(w) + \mathbf{M}_2(w), \quad \forall w \in L_w^{p_\mu}(\Omega;\delta^\alpha) \cap L^{q_1}(\Omega;\delta^{\alpha-1}), \\ &\tilde{\mathbf{M}}(w) := \tilde{\mathbf{M}}_1(w) + \mathbf{M}_2(w), \quad \forall w \in L_w^{p_\mu}(\Omega;\delta^\alpha) \cap L^{q_1}(\Omega;\delta^{\alpha-1}). \end{aligned}$$

**Step 1**: Upper bound for $g(w + \tilde\varrho\mathbb{K}_\mu[\tilde\nu] + \tilde\sigma\mathbb{G}_\mu[\tilde\tau])$ in $L^1(\Omega;\delta^\alpha)$ with $w \in L_w^{p_\mu}(\Omega;\delta^\alpha) \cap L^{q_1}(\Omega;\delta^{\alpha-1})$.

For $\lambda > 0$, set $\tilde B_\lambda := \{x \in \Omega : |w| + \tilde\varrho\mathbb{K}_\mu[|\tilde\nu|] + \tilde\sigma\mathbb{G}_\mu[|\tilde\tau|] > \lambda\}$ and $b(\lambda) := \int_{\tilde B_\lambda} \delta^\alpha dx$. We write

$$\begin{aligned} \|g(w + \tilde\varrho\mathbb{K}_\mu[\tilde\nu] + \tilde\sigma\mathbb{G}_\mu[\tilde\tau])\|_{L^1(\Omega;\delta^\alpha)} \leq &\int_{\tilde B_1} g(|w| + \tilde\varrho\mathbb{K}_\mu[|\tilde\nu|] + \tilde\sigma\mathbb{G}_\mu[|\tilde\tau|])\delta^\alpha dx \\ &+ \int_{\tilde B_1^c} g(|w| + \tilde\varrho\mathbb{K}_\mu[|\tilde\nu|] + \tilde\sigma\mathbb{G}_\mu[|\tilde\tau|])\delta^\alpha dx \\ &- \int_{\tilde B_1} g(-|w| - \tilde\varrho\mathbb{K}_\mu[|\tilde\nu|] - \tilde\sigma\mathbb{G}_\mu[|\tilde\tau|])\delta^\alpha dx \\ &- \int_{\tilde B_1^c} g(-|w| - \tilde\varrho\mathbb{K}_\mu[|\tilde\nu|] - \tilde\sigma\mathbb{G}_\mu[|\tilde\tau|])\delta^\alpha dx \\ &=: I + II + III + IV. \end{aligned} \tag{5.11}$$

We first estimate $I$. Since $g \in C(\mathbb{R}_+)$ is nondecreasing, one gets

$$I = b(1)g(1) + \int_1^\infty b(s)dg(s).$$

Since $g$ is bounded, there exists an increasing sequence of real positive number $\{\ell_j\}$ such that

$$\lim_{j\to\infty} \ell_j = \infty \quad \text{and} \quad \lim_{j\to\infty} \ell_j^{-p_\mu} g(\ell_j) = 0. \tag{5.12}$$

Observe that

$$\int_1^\infty b(s)dg(s) = \lim_{j\to\infty} \int_\lambda^{\ell_j} b(s)dg(s).$$

On the other hand, by (2.2) one gets, for every $s > 0$,

$$a(s) \leq \||w| + \tilde\varrho\mathbb{K}_\mu[|\tilde\nu|] + \tilde\sigma\mathbb{G}_\mu[|\tilde\tau|]\|^{p_\mu}_{L_w^{p_\mu}(\Omega;\delta^\alpha)} s^{-p_\mu} \leq C(\mathbf{M}_1(w) + \tilde\varrho + \tilde\sigma)^{p_\mu} s^{-p_\mu} \tag{5.13}$$

where $C = C(N,\mu,\Omega)$. Using (5.13), we obtain

$$\begin{aligned}
&b(1)g(1) + \int_1^{\ell_j} b(s)dg(s) \\
&\qquad \le C(\mathbf{M}_1(w) + \tilde{\varrho} + \tilde{\sigma})^{p_\mu} g(1) + C(\mathbf{M}_1(w) + \tilde{\varrho} + \tilde{\sigma})^{p_\mu} \int_1^{\ell_j} s^{-p_\mu} dg(s) \\
&\qquad \le C(\mathbf{M}_1(w) + \tilde{\varrho} + \tilde{\sigma})^{p_\mu} \ell_j^{-p_\mu} g(\ell_j) + Cp_\mu(\mathbf{M}_1(w) + \tilde{\varrho} + \tilde{\sigma})^{p_\mu} \int_1^{\ell_j} s^{-1-p_\mu} g(s)ds.
\end{aligned}$$

By virtue of (5.12), letting $j \to \infty$ yields

$$I \le Cp_\mu(\mathbf{M}_1(w) + \tilde{\varrho} + \tilde{\sigma})^{p_\mu} \int_1^{\infty} s^{-1-p_\mu} g(s)ds. \tag{5.14}$$

Similarly we have

$$III \le -Cp_\mu(\mathbf{M}_1(w) + \tilde{\varrho} + \tilde{\sigma})^{p_\mu} \int_1^{\infty} s^{-1-p_\mu} g(-s)ds.$$

To handle the remaining terms $II, III$, without lost of generality, we assume $q_1 \in (1, \frac{N+\alpha-1}{N+\alpha-2})$. Since $g$ satisfies condition (1.40), it follows that

$$\begin{aligned}
\max\{II, IV\} &\le a_1 \int_{\tilde{B}_1^c} (|w| + \tilde{\varrho}\mathbb{K}_\mu[|\tilde{\nu}|] + \tilde{\sigma}\mathbb{K}_\mu[|\tilde{\tau}|])^{q_1} \delta^\alpha dx + b_1 \int_{\tilde{B}_1^c} \delta^\alpha dx \\
&\le a_1 C \int_\Omega |w|^{q_1} \delta^\alpha dx + a_1 c_{34}(\tilde{\varrho}^{q_1} + \tilde{\sigma}^{q_1}) + b_1 C \\
&\le a_1 C\mathbf{M}_2(w)^{q_1} + a_1 C(\tilde{\varrho}^{q_1} + \tilde{\sigma}^{q_1}) + b_1 C
\end{aligned} \tag{5.15}$$

where $C = C(N,\mu,\Omega)$.

Combining (5.11), (5.14) and (5.15) yields

$$\|g(w + \tilde{\varrho}\mathbb{K}_\mu[\tilde{\nu}] + \tilde{\sigma}\mathbb{G}_\mu[\tilde{\tau}])\|_{L^1(\Omega;\delta^\alpha)} \le C\Lambda_g \mathbf{M}_1(w)^{p_\mu} + a_1 C\mathbf{M}_2(w)^{q_1} + b_1 C + d_{\tilde{\varrho},\tilde{\sigma}} \tag{5.16}$$

where $d_{\tilde{\varrho},\tilde{\sigma}} = C\Lambda_g(\tilde{\varrho}^{p_\mu} + \tilde{\sigma}^{p_\mu}) + a_1 C(\tilde{\varrho}^{q_1} + \tilde{\sigma}^{q_1})$.

**Step 2:** Estimates on $\mathbf{M}_1, \mathbf{M}_2$ and $\mathbf{M}$.

From (2.13), we have

$$\begin{aligned}
\tilde{\mathbf{M}}_1(\mathbb{S}(w)) &= \|\mathbb{G}_\mu[g(w + \tilde{\varrho}\mathbb{K}_\mu[\tilde{\nu}] + \tilde{\sigma}\mathbb{G}_\mu[\tilde{\tau}])\|_{L_w^{p_\mu}(\Omega;\delta^\alpha)} \\
&\le C\|g(w + \tilde{\varrho}\mathbb{K}_\mu[\tilde{\nu}] + \tilde{\sigma}\mathbb{G}_\mu[\tilde{\tau}])\|_{L^1(\Omega;\delta^\alpha)}.
\end{aligned} \tag{5.17}$$

It follows that

$$\tilde{\mathbf{M}}_1(\mathbb{S}(w)) \le C\Lambda_g \mathbf{M}_1(w)^{p_\mu} + a_1 C\mathbf{M}_2(w)^{q_1} + b_1 C + Cd_{\tilde{\varrho},\tilde{\sigma}}. \tag{5.18}$$

Applying (2.13), we get

$$\begin{aligned}
\mathbf{M}_2(\mathbb{S}(w)) &= \|\mathbb{G}_\mu[g(w + \tilde{\varrho}\mathbb{K}_\mu[\tilde{\nu}] + \tilde{\sigma}\mathbb{G}_\mu[\tilde{\tau}])\|_{L^{q_1}(\Omega;\delta^{\alpha-1})} \\
&\le C\|g(w + \tilde{\varrho}\mathbb{K}_\mu[\tilde{\nu}] + \tilde{\sigma}\mathbb{G}_\mu[\tilde{\tau}])\|_{L^1(\Omega;\delta^\alpha)},
\end{aligned}$$

which implies

$$\mathbf{M}_2(\mathbb{S}(w)) \le C\Lambda_g \mathbf{M}_1(w)^{p_\mu} + a_1 C\mathbf{M}_2(w)^{q_1} + b_1 C + Cd_{\tilde{\varrho},\tilde{\sigma}}. \tag{5.19}$$

Consequently,

$$\tilde{\mathbf{M}}(\mathbb{S}(w)) \leq C\Lambda_g \mathbf{M}_1(w)^{p_\mu} + a_1 C\mathbf{M}_2(w)^{q_1} + b_1 C + Cd_{\tilde{\varrho},\tilde{\sigma}}. \tag{5.20}$$

Similarly, we can show that

$$\mathbf{M}(\tilde{\mathbb{S}}(w)) \leq \tilde{C}\Lambda_{\tilde{g}} \tilde{\mathbf{M}}_1(w)^{p_\mu} + a_1 \tilde{C}\mathbf{M}_2(w)^{q_1} + b_1 \tilde{C} + \tilde{C}d_{\varrho,\sigma} \tag{5.21}$$

where $\tilde{C}$ is a positive constant. Define the functions $\eta$ and $\tilde{\eta}$ as follows

$$\eta(\lambda) := \max\{C\Lambda_g, \tilde{C}\Lambda_{\tilde{g}}\}\lambda^{p_\mu} + \max\{C,\tilde{C}\}a_1\lambda^{q_1} + \max\{C,\tilde{C}\}b_1 + \max\{Cd_{\tilde{\varrho},\tilde{\sigma}}, \tilde{C}d_{\varrho,\sigma}\}$$
$$\tilde{\eta}(\lambda) := \max\{C\Lambda_g, \tilde{C}\Lambda_{\tilde{g}}\}\lambda^{p_\mu} + \max\{C,\tilde{C}\}a_1\lambda^{q_1} + \max\{C,\tilde{C}\}b_1 + \max\{Cd_{\tilde{\varrho},\tilde{\sigma}}, \tilde{C}d_{\varrho,\sigma}\}$$

where $C$ and $\tilde{C}$ are the constants in (5.20) and (5.21) respectively. By (5.20) and (5.21), we deduce

$$\tilde{\mathbf{M}}(\mathbb{S}(w)) \leq \eta(\mathbf{M}(w)) \quad \text{and} \quad \mathbf{M}(\tilde{\mathbb{S}}(w)) \leq \tilde{\eta}(\tilde{\mathbf{M}}(w)).$$

Since $p_\mu > 1$ and $q_1 > 1$, there exist $\varrho_* > 0$ and $b_* > 0$ depending on $N$, $\mu$, $\Omega$, $\Lambda_g$, $\Lambda_{\tilde{g}}$, $a_1$, $q_1$ such that for any $\varrho, \tilde{\varrho} \in (0, \varrho_*)$ and $b_1 \in (0, b_*)$ there exist $\lambda_* > 0$ and $\tilde{\lambda}_* > 0$ such that

$$\eta(\lambda_*) = \tilde{\lambda}_* \quad \text{and} \quad \tilde{\eta}(\tilde{\lambda}_*) = \lambda_*.$$

Here $\lambda_*$ and $\tilde{\lambda}_*$ depend on $N$, $\mu$, $\Omega$, $\Lambda_g$, $\Lambda_{\tilde{g}}$, $a_1$, $q_1$. Therefore,

$$\begin{aligned} &\mathbf{M}(w) \leq \lambda_* \Longrightarrow \tilde{\mathbf{M}}(\mathbb{S}(w)) \leq \tilde{\lambda}_* \\ &\tilde{\mathbf{M}}(w) \leq \tilde{\lambda}_* \Longrightarrow \mathbf{M}(\tilde{\mathbb{S}}(w)) \leq \lambda_*. \end{aligned} \tag{5.22}$$

**Step 3:** To apply Schauder fixed point theorem.

For $w_1, w_2 \in L^1(\Omega)$, put

$$\mathbb{T}(w_1, w_2) := (\mathbb{S}(w_2), \tilde{\mathbb{S}}(w_1)), \tag{5.23}$$

$$\mathcal{D} := \{(\varphi, \tilde{\varphi}) \in L^1_+(\Omega) \times L^1_+(\Omega) : \mathbf{M}(\varphi) \leq \lambda_* \text{ and } \tilde{\mathbf{M}}(\tilde{\varphi}) \leq \tilde{\lambda}_*\}.$$

Clearly, $\mathcal{D}$ is a convex subset of $L^1(\Omega) \times L^1(\Omega)$. We shall show that $\mathcal{D}$ is a closed subset of $(L^1(\Omega))^2$. Indeed, let $\{(\varphi_m, \tilde{\varphi}_m)\}$ be a sequence in $\mathcal{D}$ converging to $(\varphi, \tilde{\varphi})$ in $(L^1(\Omega))^2$. Obviously, $\varphi \geq 0$ and $\tilde{\varphi} \geq 0$. We can extract a subsequence, still denoted by the same notation, such that $(\varphi_m, \tilde{\varphi}_m) \to (\varphi, \tilde{\varphi})$ a.e. in $\Omega$. Consequently, by Fatou's lemma,

$$\mathbf{M}_i(\varphi) \leq \liminf_{m\to\infty} \mathbf{M}_i(\varphi_m), \quad \mathbf{M}_i(\tilde{\varphi}) \leq \liminf_{m\to\infty} \mathbf{M}_i(\tilde{\varphi}_m)$$

for $i = 1, 2$. It follows that $\mathbf{M}(\varphi) \leq \lambda_*$ and $\mathbf{M}(\tilde{\varphi}) \leq \tilde{\lambda}_*$. So $(\varphi, \tilde{\varphi}) \in \mathcal{D}$ and therefore $\mathcal{D}$ is a closed subset of $L^1(\Omega) \times L^1(\Omega)$.

Clearly, $\mathbb{T}$ is well defined in $\mathcal{D}$. For $(w, \tilde{w}) \in \mathcal{D}$, we get $\mathbf{M}(w) \leq \lambda_*$ and $\mathbf{M}(\tilde{w}) \leq \tilde{\lambda}_*$, hence $\tilde{\mathbf{M}}(\mathbb{S}(w)) \leq \tilde{\lambda}_*$ and $\mathbf{M}(\mathbb{S}(\tilde{w})) \leq \lambda_*$. It follows that $\mathbb{T}(\mathcal{D}) \subset \mathcal{D}$.

We observe that $\mathbb{T}$ is continuous. Indeed, if $w_m \to w$ and $\tilde{w}_m \to \tilde{w}$ as $m \to \infty$ in $L^1(\Omega)$ then since $g, \tilde{g} \in C(\mathbb{R}) \cap L^\infty(\mathbb{R})$, it follows that

$$g(\tilde{w}_m + \tilde{\varrho}\mathbb{K}_\mu[\nu] + \tilde{\sigma}\mathbb{G}_\mu[\tilde{\tau}]) \to g(\tilde{w} + \tilde{\varrho}\mathbb{K}_\mu[\nu] + \tilde{\sigma}\mathbb{G}_\mu[\tilde{\tau}]) \quad \text{in } L^1(\Omega; \delta^\alpha),$$

and

$$\tilde{g}(w_m + \varrho\mathbb{K}_\mu[\nu] + \sigma\mathbb{G}_\mu[\tau]) \to \tilde{g}(w + \varrho\mathbb{K}_\mu[\nu] + \sigma\mathbb{G}_\mu[\tau]) \quad \text{in } L^1(\Omega; \delta^\alpha)$$

as $m \to \infty$. By (2.13), $\mathbb{S}(\tilde{w}_m) \to \mathbb{S}(\tilde{w})$ and $\tilde{\mathbb{S}}(w_m) \to \tilde{\mathbb{S}}(w)$ as $m \to \infty$ in $L^1(\Omega)$. Thus $\mathbb{T}(w_m, \tilde{w}_m) \to \mathbb{T}(w, \tilde{w})$ in $L^1(\Omega) \times L^1(\Omega)$.

We next show that $\mathbb{T}$ is a compact operator. Let $\{(w_m, \tilde{w}_n)\} \subset \mathcal{D}$ and for each $m \geq 1$, put $\psi_m = \mathbb{S}(\tilde{w}_m)$ and $\tilde{\psi}_m = \tilde{\mathbb{S}}(w_m)$. Hence $\{\Delta\psi_m\}$ and $\{\Delta\tilde{\psi}_m\}$ are uniformly bounded in $L^p(G)$ for every subset $G \Subset \Omega$. Therefore $\{\psi_m\}$ is uniformly bounded in $W^{1,p}(G)$. Consequently, there exists a subsequence, still denoted by the same notation, and functions $\psi, \tilde{\psi}$ such that $(\psi_m, \tilde{\psi}_m) \to (\psi, \tilde{\psi})$ a.e. in $\Omega$. By dominated convergence theorem, $(\psi_m, \tilde{\psi}_m) \to (\psi, \tilde{\psi})$ in $L^1(\Omega) \times L^1(\Omega)$. Thus $\mathbb{T}$ is compact.

By Schauder fixed point theorem there is $(u, v) \in \mathcal{D}$ such that $\mathbb{T}(u, v) = (u, v)$. □

**Proof of Theorem H.**I. Let $\{g_n\}$ and $\{\tilde{g}_n\}$ be the sequences of continuous, nondecreasing functions defined on $\mathbb{R}$ such that

$$\begin{aligned} &g_n(0) = g(0),\ |g_n| \leq |g_{n+1}| \leq |g|,\ \sup_{\mathbb{R}} |g_n| = n \text{ and } \lim_{n\to\infty} \|g_n - g\|_{L^\infty_{\text{loc}}(\mathbb{R})} = 0, \\ &\tilde{g}_n(0) = \tilde{g}(0),\ |\tilde{g}_n| \leq |\tilde{g}_{n+1}| \leq |\tilde{g}|,\ \sup_{\mathbb{R}} |\tilde{g}_n| = n \text{ and } \lim_{n\to\infty} \|\tilde{g}_n - \tilde{g}\|_{L^\infty_{\text{loc}}(\mathbb{R})} = 0. \end{aligned} \tag{5.24}$$

Due to Lemma 5.1, there exist $\lambda_*$, $\tilde{\lambda}_*$, $b_* > 0$ and $\varrho_* > 0$ depending on $N$, $\mu$, $\Omega$, $\Lambda_g$, $\Lambda_{\tilde{g}}$, $a_1$, $q_1$ such that for every $b_1 \in (0, b_*)$, $\tilde{\varrho}, \varrho \in (0, \varrho_*)$ and $n \geq 1$ there exists a solution $(w_n, \tilde{w}_n) \in \mathcal{D}$ of

$$\begin{cases} -L_\mu w_n = g_n(\tilde{w}_n + \tilde{\varrho}\mathbb{K}_\mu[\tilde{\nu}] + \tilde{\sigma}\mathbb{G}_\mu[\tilde{\tau}]) & \text{in } \Omega, \\ -L_\mu \tilde{w}_n = \tilde{g}_n(w_n + \varrho\mathbb{K}_\mu[\nu] + \sigma\mathbb{G}_\mu[\tau]) & \text{in } \Omega, \\ \operatorname{tr}(w_n) = \operatorname{tr}(\tilde{w}_n) = 0. \end{cases} \tag{5.25}$$

For each $n$, set $u_n = w_n + \varrho\mathbb{K}_\mu[\nu] + \sigma\mathbb{G}_\mu[\tau]$ and $v_n = \tilde{w}_n + \tilde{\varrho}\mathbb{K}_\mu[\tilde{\nu}] + \tilde{\sigma}\mathbb{G}_\mu[\tilde{\tau}]$. Then

$$-\int_\Omega u_n L_\mu \phi dx = \int_\Omega g_n(v_n)\phi dx + \sigma \int_\Omega \phi d\tau - \varrho \int_\Omega \mathbb{K}_\mu[\nu] L_\mu \phi \, dx \quad \forall \phi \in \mathbf{X}_\mu(\Omega), \tag{5.26}$$

$$-\int_\Omega v_n L_\mu \phi dx = \int_\Omega \tilde{g}_n(u_n)\phi dx + \tilde{\sigma} \int_\Omega \phi d\tilde{\tau} - \tilde{\varrho} \int_\Omega \mathbb{K}_\mu[\tilde{\nu}] L_\mu \phi \, dx \quad \forall \phi \in \mathbf{X}_\mu(\Omega). \tag{5.27}$$

Since $\{(w_n, \tilde{w}_n)\} \subset \mathcal{D}$ and the fact that $\Lambda_{g_n} \leq \Lambda_g$, we obtain from (5.16) that

$$\|g_n(v_n)\|_{L^1(\Omega;\delta^\alpha)} \leq C\Lambda_g \lambda_*^{p_\mu} + a_1 C \lambda_*^{q_1} + b_* C + d_{\varrho_*} \tag{5.28}$$

Hence the sequence $\{g(v_n)\}$ is uniformly bounded in $L^1(\Omega; \delta^\alpha)$. Since $\{(w_n, \tilde{w}_n)\} \subset \mathcal{D}$, the sequence $\{\frac{\mu}{\delta^2} w_n\}$ and $\{\frac{\mu}{\delta^2} \tilde{w}_n\}$ are uniformly bounded in $L^{q_1}(G)$ for every subset $G \Subset \Omega$. As a consequence, $\{\Delta w_n\}$ and $\{\Delta \tilde{w}_n\}$ are uniformly bounded in $L^1(G)$ for every subset $G \Subset \Omega$. By regularity results for elliptic equations, there exist subsequences, still denoted by the same notations, and functions $w$ and $\tilde{w}$ such that $(w_n, \tilde{w}_n) \to (w, \tilde{w})$ a.e. in $\Omega$. Therefore $(u_n, v_n) \to (u, v)$ a.e. in $\Omega$ with $u = w + \varrho\mathbb{K}_\mu[\nu] + \sigma\mathbb{G}_\mu[\tau]$ and $u = \tilde{w} + \tilde{\varrho}\mathbb{K}_\mu[\tilde{\nu}] + \tilde{\sigma}\mathbb{G}_\mu[\tilde{\tau}]$. Moreover $(\tilde{g}_n(u_n), g_n(v_n)) \to (\tilde{g}(u), g(v))$ a.e. in $\Omega$.

We show that $u_n \to u$ in $L^1(\Omega; \delta^\alpha)$. Since $\{w_n\}$ is uniformly bounded in $L^{q_1}(\Omega; \delta^{\alpha-1})$, by (2.14), we derive that $\{u_n\}$ is uniformly bounded in $L^{q_1}(\Omega; \delta^\alpha)$. Due to Holder inequality, $\{u_n\}$ is uniformly integrable with respect to $\delta^\alpha dx$. We invoke Vitali convergence theorem to derive that $u_n \to u$ in $L^1(\Omega; \delta^\alpha)$. Similarly, one can prove that $v_n \to v$ in $L^1(\Omega; \delta^\alpha)$.

We next prove that $g_n(v_n) \to g(v)$ in $L^1(\Omega;\delta^\alpha)$. For $\lambda>0$ and $n\in\mathbb{N}$ set $B_{n,\lambda} := \{x\in\Omega : |v_n|>\lambda\}$ and $b_n(\lambda) := \int_{B_{n,\lambda}} \delta^\alpha dx$. For any Borel set $E\subset\Omega$,

$$\begin{aligned}\int_E g_n(v_n)\delta^\alpha dx &= \int_{E\cap B_{n,\lambda}} g_n(v_n)\delta^\alpha dx + \int_{E\cap B_{n,\lambda}^c} g_n(v_n)\delta^\alpha dx \\ &\le \int_{B_{n,\lambda}} g_n(v_n)\delta^\alpha dx + m_{g,\lambda}\int_E \delta^\alpha dx \\ &\le b_n(\lambda)g_n(\lambda) + \int_\lambda^\infty b_n(s)dg_n(s) + m_{g,\lambda}\int_E \delta^\alpha dx.\end{aligned} \tag{5.29}$$

where $m_{g,\lambda} := \sup_{[0,\lambda]} g$. By proceeding as in the proof of Lemma 5.1 in order to get (5.14), we deduce

$$b_n(\lambda)g_n(\lambda) + \int_\lambda^\infty b_n(s)dg_n(s) \le C\int_\lambda^\infty s^{-1-p_\mu} g_n(s)ds \le C\int_\lambda^\infty s^{-1-p_\mu} g(s)ds \tag{5.30}$$

where $C$ depends on $N$, $\mu$, $\Lambda_g$, $\Lambda_{\tilde g}$, $a_1$, $q_1$. Note that the term on the right hand-side of (5.30) tends to 0 as $\lambda\to\infty$. Take arbitrarily $\varepsilon>0$, there exists $\lambda>0$ such that the right hand-side of (5.30) is smaller than $\frac{\varepsilon}{2}$. Fix such $\lambda$ and put $\eta = \frac{\varepsilon}{2m_{g,\lambda}}$. Then, by (5.29),

$$\int_E \delta(x)^\alpha dx \le \eta \Longrightarrow \int_E g_n(v_n)\delta(x)^\alpha dx < \varepsilon.$$

Therefore the sequence $\{g_n(v_n)\}$ is uniformly integrable with respect to $\delta^\alpha dx$. Due to Vitali convergence theorem, we deduce that $g_n(v_n)\to g(v)$ in $L^1(\Omega;\delta^\alpha)$.

By sending $n\to\infty$ in each term of (5.26) we obtain

$$-\int_\Omega u L_\mu\phi dx = \int_\Omega g(v)\phi dx + \sigma\int_\Omega \phi\, d\tau - \varrho\int_\Omega \mathbb{K}_\mu[\nu]L_\mu\phi\, dx, \quad \forall\phi\in\mathbf{X}_\mu(\Omega). \tag{5.31}$$

Similarly, one can show that $\tilde g_n(u_n)\to\tilde g(u)$ in $L^1(\Omega;\delta^\alpha)$. By letting $n\to\infty$ in (5.27), we get

$$-\int_\Omega v L_\mu\phi dx = \int_\Omega \tilde g(u)\phi dx + \tilde\sigma\int_\Omega \phi\, d\tilde\tau - \tilde\varrho\int_\Omega \mathbb{K}_\mu[\tilde\nu]L_\mu\phi\, dx, \quad \forall\phi\in\mathbf{X}_\mu(\Omega). \tag{5.32}$$

Thus $(u,v)$ is a solution of (1.27). □

5.3. **Source case : sublinearity.** We next deal with the case where $g$ and $\tilde g$ are sublinear.

**Proof of Theorem H.**II. The proof is similar to that of Lemma 5.1, also based on Schauder fixed point theorem. So we point out only the main modifications. Let $\mathbb{S}$ and $\tilde{\mathbb{S}}$ be the operators defined in (5.10). Put

$$\mathbf{N}_1(w) := \|w\|_{L^{q_1}(\Omega;\delta^{\alpha-1})}, \quad \forall w\in L^{q_1}(\Omega;\delta^{\alpha-1}),$$

$$\mathbf{N}_2(w) := \|w\|_{L^{1}(\Omega;\delta^{\alpha-1})}, \quad \forall w\in L^{1}(\Omega;\delta^{\alpha-1}).$$

Combining (2.13), (2.14) and (1.42) leads to

$$\mathbf{N}_2(\mathbb{S}(w)) \le a_2 C\mathbf{N}_1(w)^{q_1} + C(\tilde\varrho^{q_1} + \tilde\sigma^{q_1} + b_2).$$

On the other hand

$$\mathbf{N}_1(\tilde{\mathbb{S}}(w)) \le a_2 C\mathbf{N}_2(w)^{q_2} + C(\varrho^{q_2} + \sigma^{q_2} + b_2)..$$

Define

$$\xi_1(\lambda) := a_2 C\lambda^{q_1} + C(\tilde{\varrho}^{q_1} + \tilde{\sigma}^{q_1} + b_2),$$
$$\xi_2(\lambda) := a_2 C\lambda^{q_2} + C(\varrho^{q_2} + \sigma^{q_2} + b_2).$$

Then

$$\mathbf{N}_2(\mathbb{S}(w)) \le \xi_1(\mathbf{N}_1(w)) \quad \text{and} \quad \mathbf{N}_1(\tilde{\mathbb{S}}(w)) \le \xi_2(\mathbf{N}_2(w)),$$

If $q_1 q_2 < 1$ then we can find $\lambda_1$ and $\lambda_2$ such that $\xi_1(\lambda_1) = \lambda_2$ and $\xi_2(\lambda_2) = \lambda_1$. Thus if $\mathbf{N}_1(w) < \lambda_1$ then $\mathbf{N}_2(\mathbb{S}(w)) < \lambda_2$ and if $\mathbf{N}_2(w) < \lambda_2$ then $\mathbf{N}_1(\tilde{\mathbb{S}}(w)) < \lambda_1$.

If $q_1 q_2 = 1$ and $a_2$ small enough we can find, $\lambda_1$ and $\lambda_2$ such that $\xi_1(\lambda_1) = \lambda_2$ and $\xi_2(\lambda_2) = \lambda_1$.

The rest of the proof can be proceeded as in the proof of Lemma 5.1 and the proof of Theorem H.I. and we omit it. □

5.4. **Source case : subcriticality and sublinearity.**

**Proof of Theorem H.**III.

Set

$$\mathbf{N}(w) := \|w\|_{L^{q_1}(\Omega;\delta^{\alpha-1})}, \quad \forall w \in L^{q_1}(\Omega;\delta^{\alpha-1}).$$

By an argument similar to the proof of Lemma 5.1 and and Theorem H.II, we get

$$\mathbf{N}(\mathbb{S}(w)) \le C\Lambda_g \mathbf{M}_1(w)^{p_\mu} + a_1 C\mathbf{M}_2(w)^{q_1} + b_1 C + d_{\tilde{\varrho},\tilde{\sigma}}.$$

On the other hand

$$\mathbf{M}(\tilde{\mathbb{S}}(w)) \le a_2 C\mathbf{N}(w)^{q_2} + C(\varrho^{q_2} + \sigma^{q_2} + b_2).$$

Set

$$\hat{\xi}_1(\lambda) := C\Lambda_g \lambda^{p_\mu} + a_1 C\lambda^{q_1} + b_1 C + d_{\tilde{\varrho},\tilde{\sigma}},$$
$$\hat{\xi}_2(\lambda) := a_2 C\lambda^{q_2} + C(\varrho^{q_2} + \sigma^{q_2} + b_2).$$

Then

$$\mathbf{N}(\mathbb{S}(w)) \le \hat{\xi}_1(\mathbf{M}(w)) \quad \text{and} \quad \mathbf{M}(\tilde{\mathbb{S}}(w)) \le \hat{\xi}_2(\mathbf{N}(w)).$$

We consider there cases.

**Case 1:** $q_1 q_2 > 1$. Since $p_\mu > q_1$, it follows that $p_\mu q_2 > 1$. Therefore there exist $b_* > 0$ and $\varrho_* > 0$ such that for $b_1, b_2 \in (0, b_*)$ and $\varrho \in (0, \varrho_*)$ one can find $\lambda_1 > 0$ and $\lambda_2 > 0$ satisfying

$$\hat{\xi}_1(\lambda_1) = \lambda_2 \quad \text{and} \quad \hat{\xi}_2(\lambda_2) = \lambda_1. \tag{5.33}$$

**Case 2:** $p_\mu q_2 = 1$. In this case, there exist $a_* > 0$ such that if $a_2 \in (0, a_*)$ then for every $\varrho > 0$ and $\tilde{\varrho} > 0$ one can find $\lambda_1 > 0$ and $\lambda_2 > 0$ satisfying (5.33).

**Case 3:** $p_\mu q_2 < 1$. In this case for every $\varrho > 0$ and $\tilde{\varrho} > 0$ one can find $\lambda_1 > 0$ and $\lambda_2 > 0$ such that (5.33) holds.

Hence, in any case,

$$\mathbf{M}(w) \le \lambda_1 \Longrightarrow \mathbf{N}(\mathbb{S}(w)) \le \hat{\xi}_1(\lambda_1) = \lambda_2$$
$$\mathbf{N}(w) \le \lambda_2 \Longrightarrow \mathbf{M}(\tilde{\mathbb{S}}(w)) \le \hat{\xi}_2(\lambda_2) = \lambda_1.$$

The rest of the proof can be proceeded as in the proof of Lemma 5.1 and the proof of Theorem H.II. and we omit it. □

## References

[1] D. R. Adams, L.I. Heberg, *Function Spaces and Potential Theory,* Grundlehren der Mathematischen Wisenschaften 31, Springer-Verlag (1999).

[2] A. Ancona, *Negatively curved manifolds, elliptic operators and the Martin boundary*, Ann. of Math. (2), **125** (1987), 495-536.

[3] C. Bandle, V. Moroz and W. Reichel, *Boundary blowup type sub-solutions to semilinear elliptic equations with Hardy potential*, J. London Math. Soc. **2** (2008), 503-523.

[4] M. F. Bidaut-Véron and L. Vivier, *An elliptic semilinear equation with source term involving boundary measures: the subcritical case*, Rev. Mat. Iberoamericana **16** (2000), 477-513.

[5] M. F. Bidaut-Véron, G. Hoang, Q. H. Nguyen, L. Véron, *An elliptic semilinear equation with source term and boundary measure data: the supercritical case*, Journal of Functional Analysis **269** (2015), 1995-2017.

[6] M.F. Bidaut-Véron and C. Yarur, *Semilinear elliptic equations and systems with measure data: existence and a priori estimates*, Adv. Differential Equations **7** (2002), 257-296

[7] S. Filippas, L. Moschini and A. Tertikas, *Sharp two-sided heat kernel estimates for critical Schrodinger operators on bounded domains*, Commun. Math. Phys. **273** (2007), 237-281.

[8] M. García-Huidobro, R. Manásevich, E. Mitidieri, Enzo, C. S. Yarur, *Existence and nonexistence of positive singular solutions for a class of semilinear elliptic systems*, Arch. Rational Mech. Anal. **140** (1997), no. 3, 253-284.

[9] K. T. Gkikas and L. Véron, *Boundary singularities of solutions of semilinear elliptic equations with critical Hardy potentials*, Nonlinear Anal. **121** (2015), 469-540.

[10] A. Gmira and L. Véron, *Boundary singularities of solutions of some nonlinear elliptic equations*, Duke Math. J. **64** (1991), 271-324.

[11] N.J. Kalton and I.E. Verbitsky, *Nonlinear equations and weighted norm inequality*, Trans. Amer. Math. Soc. **351** (1999) 3441-3497.

[12] M. Marcus, V. J. Mizel and Y. Pinchover, *On the best constant for Hardy's inequality in* $\mathbb{R}^N$, Trans. Amer. Math. Soc. **350** (1998), 3237-3255.

[13] M. Marcus and V. Moroz, *Moderate solutions of semilinear elliptic equations with Hardy potential under minimal restrictions on the potential*, http://arxiv.org/abs/1603.09265.

[14] M. Marcus and P. T. Nguyen, *Moderate solutions of semilinear elliptic equations with Hardy potential*, Ann. Inst. H. Poincaré Anal. Non Linéaire **34** (2017), 69–88.

[15] M. Marcus and L. Véron, *The boundary trace of positive solutions of semilinear elliptic equations: the subcritical case*, Arch. Rational Mech. Anal. **144** (1998), 201-231.

[16] M. Marcus and L. Véron, *The boundary trace of positive solutions of semilinear elliptic equations: the supercritical case*, J. Math. Pures Appl. (9) **77** (1998), 481-524.

[17] M. Marcus and L. Véron, *Nonlinear second order elliptic equations involving measures*, De Gruyter Series in Nonlinear Analysis and Applications, 2013.

[18] P.-T. Nguyen, *Elliptic equations with Hardy potential and subcritical source term*, Calc. Var. Partial Differential Equations **56** (2017), no. 2, Art. 44, 28 pp.

[19] P. Quittner and P. Souplet, *Superlinear parabolic problems. Blow-up, global existence and steady states*, Birkhäuser Verlag, Basel, 2007. xii+584 pp.

Konstantinos T. Gkikas, Centro de Modelamiento Matemático (UMI 2807 CNRS), Universidad de Chile, Casilla 170 Correo 3, Santiago, Chile.

Phuoc-Tai Nguyen, Departamento de Matemática, Pontificia Universidad Católica de Chile, Santiago, Chile and Department of Mathematics and Statistics, Masaryk University, Brno, Czech Republic